\documentclass[12pt]{amsart}
\usepackage{epsfig,amssymb,epic,eepic,color}
\usepackage[all]{xy}
\numberwithin{equation}{section}
 \usepackage[T1]{fontenc}        
\newtheorem{thm}{Theorem}[section]

\newtheorem{prop}[thm]{Proposition}

\newtheorem{lemme}[thm]{Lemma}
\newtheorem{cor}[thm]{Corollary}

\newtheorem{remarque}[thm]{Remark}
\newtheorem{remarques}[thm]{Remarks}

\newtheorem{rien}[thm]{}

\numberwithin{equation}{section}
\newcommand{\be}{\begin{enumerate}}
\newcommand{\ee}{\end{enumerate}}
\newcommand{\bi}{\begin{itemize}}
\newcommand{\ei}{\end{itemize}}

\def\R{\mathbb{R}}

\def\D{\mathbb{D}}
\def\F{\mathcal{F}}
\def\G{\mathcal{G}}
\def\XX{\mathcal{X}}
\def\PP{\mathcal{P}}
\def\S{\mathbb{S}}

\def\om{\omega}
\def\Om{\Omega}

\def\ga{\gamma}    
\def\Ga{\Gamma}

\def\al{\alpha}
\def\be{\beta}
 
\def\de{\delta}
\def\De{\Delta}
\def\vp{\varphi}

\def\la{\lambda}
\def\La{\Lambda}

\def\si{\sigma}
\def\Si{\Sigma}

\def\ep{\varepsilon}

\def\nd{\noindent}
\def\bull{\hfill$\Box$\\}
\def\proof{\nd {\bf Proof.\ }}

\begin{document}

\vskip 1cm
\begin{center}{\sc Haefliger structures and  symplectic/contact structures
\vskip 1cm

 Fran\c cois Laudenbach \& Ga\"el Meigniez}

\end{center}
\title{}

\author{ }
\address{Universit\'e de Nantes, LMJL, UMR 6629 du CNRS, 44322 Nantes, France}
\email{francois.laudenbach@univ-nantes.fr}

\address{Universit\'e de Bretagne Sud,  LMBA, UMR 6205 du CNRS, 
BP 573, F-56017 Vannes, France}
\email{Gael.Meigniez@univ-ubs.fr}



\thanks{FL is supported by ERC Geodycon}

\begin{abstract}  For some geometries including symplectic and contact structures on an $n$-dimensional
manifold,
we introduce a two-step approach to Gromov's $h$-principle. 
From  formal geometric data, the first step builds a transversely geometric Haefliger structure of 
codimension $n$.
This step works on all manifolds, even closed. The second step, which works only on open manifolds and 
for all geometries, regularizes the intermediate Haefliger structure and produces a genuine geometric 
structure. Both steps admit relative parametric versions. The proofs borrow ideas from W. Thurston,
like jiggling and inflation. Actually, we are using a more primitive jiggling due to R. Thom.

\end{abstract}

\maketitle
\thispagestyle{empty}
\vskip 1cm

\section{Introduction}\label{intro}

We consider geometric structures on manifolds, such as the following:
symplectic structure,  contact structure,  foliation of prescribed codimension,
 immersion or  submersion to another manifold.
We recall that, in order to provide  
a given manifold $M$
 with such a structure,
Gromov's $h$-principle
consists of starting from a \emph{formal}
version of the structure on $M$ (this means a {\it non-holonomic} -- 
 that is,  non-integrable -- section of some jet space) and deforming it until it becomes \emph{genuine}
(holonomic) \cite{gromov}. In the present paper, we introduce a two-step
 approach to the
$h$-principle for such structures.

From the formal data, the first step builds
 a Haefliger structure of codimension zero on $M$, transversely geometric; this concept 
  will be explained  below.
For each of the geometries  above-mentioned,
 the first step  
 { works} for every manifold $M$, even closed.

  The second step,  which is works 
  for open manifolds only, regularizes the intermediate Haefliger structure, providing a genuine geometric structure. 
 Both steps admit 
relative  parametric versions.

An essential tool in both steps 
consists of {\it jiggling.}
We recall that Thurston's work on foliations used his famous Jiggling Lemma \cite{thu74}.
 As A. Haefliger told us 
\cite{haefliger2},
Thurston himself was aware that this lemma applies for getting some $h$-principles
 in the sense
of Gromov. 
In a not very popular paper by R. Thom \cite{thom59}, we discovered 
 a more primitive jiggling lemma
that   is  
remarkably suitable for the needs of our approach.

\medbreak 

\begin{rien}{\bf Groupoids and geometries.}\label{list-geo} {\rm According to 
 O. Veblen and J.H.C. Whitehead
 \cite{veblen_31}, a geometry in dimension $n$ is defined by an
   $n$-dimensional model manifold $X$ (often $\R^n$) and by an open  subgroupoid $\Gamma$ 
    in the groupoid 
  $\Gamma(X)$ 
{   of the germs of local $C^\infty$-diffeomorphisms of $X$;}
 here the topology on $\Gamma(X)$ is meant to be the
  {\it sheaf } topology. 
In what follows, we use the classical  notation $\Ga_n:= \Ga(\R^n)$. Here are examples of such open
subgroupoids.}
\end{rien}
\begin{enumerate}
\item When $n$ is even,  $\Gamma^{ \rm symp}_n\subset \Ga_n$ denotes the subgroupoid of germs
preserving the standard symplectic form of $\R^n$.
\item When $n$ is odd, $\Gamma^{\rm cont}_n\subset \Ga_n$ denotes the subgroupoid of germs 
preserving the standard (positive) contact structure of $\R^n$.  
\item  For $n=p+q$, one has the subgroupoid
 $\Gamma^{\rm fol}_{n,q}\subset\Gamma_n$
  preserving the standard foliation
 of codimension $q$ (whose leaves are  the $p$-planes parallel to $\R^p$).
\item  When $Y$ is any $q$-dimensional manifold and  $X=\R^p\times Y$,  
 one has the subgroupoid
 $\Gamma_n^Y\subset\Gamma(X)$ of the germs of the form
 $(x,y)\mapsto
 (f(x,y),y)$.
\end{enumerate}

\begin{rien} {\bf $\Gamma$-foliations.}
\label{gamma_fol} {\rm The next concept goes back to A. Haefliger  \cite{haefliger}.  For an open subgroupoid $\Ga$ of  
 $\Ga(X)$, a \emph{$\Gamma$-foliation} on a manifold $E$ is meant to be  a codimension-$n$ foliation 
 on $E$ equipped with a {\it transverse geometry} associated with $\Ga$ and invariant by holonomy. More precisely,
 this foliation is 
defined by a maximal atlas of submersions $(f_i:U_i\to X)$ from  open subsets
of $E$ into $X$ such that,  for every $i, j$ and every $x\in U_i\cap U_j$, there is
a germ  $\gamma_{ij}\in\Gamma$ at point  $f_j(x)$ verifying 
$$[f_i]_x=\gamma_{ij}\, [f_j]_x
$$
where [-]$_x$ stands for the germ at $x$. When $n=\dim E$, one also speaks of  
a \emph{$\Gamma$-geometry} on $E$.}
 \end{rien}
 
 Here are examples related to the previous list of groupoids. 
 The first two have been already considered by D. McDuff in \cite{mcduff}.
 \begin{enumerate}
 \item A $\Gamma^{\rm symp}_n$-foliation on $E$
 amounts to a closed differential $2$-form $\Omega$ on $E$, whose kernel
 is of codimension $n$ at every point. The closedness of $\Om$ is equivalent to the conjunction of the next  facts:
 \begin{itemize}
 \item  the  codimension-$n$ plane field $(x\in E\mapsto \ker \Om_x)$ is integrable,
 \item  $\Om$ is basic\footnote{We recall that a form $\al$ is said to be basic with respect to a foliation $
 \mathcal F$ if  the Lie derivative $L_X\al$ vanishes for every vector field 
 $X$ tangent to $\mathcal F$; this  is the infinitesimal version of the
{\it invariance by holonomy}. } with respect to that foliation,
 \item  $\Om$ is closed on a total transversal.
 \end{itemize}
   
 \item   A $\Gamma_n^{\rm cont}$-foliation
 on $E$ amounts to a  codimension-one plane field  $P$ on $E$ defined by
 an equation $A=0$ where $A$ is a differential form of degree 1, unique up to multiplying by 
  a positive\footnote{Here, we limit ourselves to co-orientable contact structures.} function, which satisfies 
    the next conditions:
 \begin{itemize}
 \item the $n$-form $A\wedge(dA)^{(n-1)/2}$ is closed and has  a codimension-$n$ kernel $K_x$ at every point $x\in E$; in particular, the field $(x\mapsto K_x)$ is integrable, tangent to a codimension-$n$ foliation denoted by $\mathcal K$;
 \item $K_x$ is a vector sub-space of $P_x$ for every $x$;
 \item $P$ is invariant by the holonomy of $\mathcal K$.
 \end{itemize}
  \item  A $\Gamma^{\rm fol}_{n,q}$-foliation on $E$
 consists of  a flag $\F\subset\G$ of
  two nested foliations of respective codimensions $n$ and $q$ in $E$ with $n>q$.
 \item A $\Gamma_n^Y$-foliation on $E$ consists of a  codimension-$n$ foliation
  and a submersion {   $w:E\to Y$} which is
 constant on every leaf (\cite{haefliger} p. 145).\hfill${}$
   \end{enumerate}
\medbreak

\begin{rien} {\bf Haefliger's $\Gamma$-structures}. \label{gamma_struct}
{\rm A. Haefliger defined a \emph{$\Gamma$-structure} as some class of cocycles   valued
 in $\Gamma$  (see \cite{haefliger}  p. 137). 
This definition, that
 makes sense on every topological space and for every topological groupoid $\Gamma$,
  allowed him
to build a classifying space $B\Gamma$ for these structures (\cite{haefliger} p. 140).
A second description (\cite{haefliger70} p. 188) is 
more suitable for our purpose
when the topological space is a manifold $M$ and when the groupoid is
an open subgroupoid $\Gamma$  in the
 groupoid of germs $\Gamma(X)$ of a $n$-manifold $X$. Here it is.}
 \end{rien}
 {\it A $\Gamma$-structure 
 on $M$ consists of a pair $\xi= (\nu,\F )$ where
\begin{itemize}
\item $\nu$ is a real vector bundle over $M$ of rank $n$, called the \emph{normal bundle};
 its total space
is denoted by   $E(\nu)$; and $Z:M\to E(\nu)$ denotes the zero section;
 \item $\F$ is a germ along $Z(M)$ of $\Gamma$-foliation 
in $E(\nu)$ transverse to every fibre of $\nu$.
\end{itemize}
}
An important feature of $\Ga$-structures is that { pulling back} by smooth maps
 (in our restricted setting)  { is} 
allowed without assuming any transversality: if $f: N\to M$ is a smooth map and $\xi$
 is a $\Ga$-structure 
on $M$,
one defines $f^*\xi$ as
 the $\Ga$-structure on $N$ whose normal bundle is $f^*\nu$ equipped with the 
$\Ga$-foliation $F^{-1}(\F)$, where $F$ is {   a} 
bundle morphism  over $f$ which is a 
fibre-to-fibre linear isomorphism.

Let $H^1(M;\Ga)$ (resp. $H^1_\nu(M;\Gamma)$)
 denote the space of the $\Gamma$-structures on $M$
(resp. {those} whose normal bundle is $\nu$). 
 It is a topological space since $\Ga$ is a topological groupoid; their elements
are denoted by $\xi=(\nu,\F)$.
In what follows, we are mainly 
interested in the case where $\dim M=n$ and $\nu$ is isomorphic to the
 tangent space $\tau M$; in that case, the elements are just denoted $\F$.
 
 In what follows, a $\Ga$-structure  on $M$ whose normal bundle 
is the tangent bundle $\tau M$
will be called a \emph{tangential} $\Ga$-structure. {    For short, when it is not ambiguous,
 we write $\mathcal F$
for $(\tau M, \mathcal F)$.  We introduced the general definition
of $\Ga$-structure  --
at least in the smooth case -- since we are going to refer to in several places.
}
\medbreak

\begin{rien} {\bf Underlying formal geometries.}\label{formal} {\rm 
 In the cases of the geometries (1), (2) and (4) given above, that is,
 for $\Gamma=\Gamma_n^{\rm symp}$, $\Gamma_n^{\rm cont}$
 or $\Gamma_n^Y$,
 every $\Gamma$-structure has an 
underlying 
 formal $\Gamma$-geometry in the sense of Gromov.  
 But, we do not intend to enter Gromov's generality. We just describe what they are.
 In the case of the geometry (3), every $\Gamma^{\rm fol}_{n,q}$-structure has an
 underlying object, somewhat formal, but more complicated than in the cases
 (1), (2), (4).

}
\end{rien}
\begin{enumerate}
\item Assume $n$ is even.
Given $\F\in H^1_{\tau M}(M;\Gamma_n^{\rm symp})$, one has
 an
associated basic closed $2$-form $\Omega$ on a neighborhood of $Z(M)$
in the total space $TM$.  Its kernel is everywhere transverse to the fibres. Therefore,  
 $\Omega$ defines
 a non-singular $2$-form
 $\omega$ on $M$ by the formula  $\omega_x:=\Omega_{Z(x)}\vert T_xM$ for every $x\in M$.
 { This is the underlying {\it formal symplectic structure}.}
 \item Assume $n$ is odd.
Given $\F\in H^1_{\tau M}(M;\Gamma_n^{\rm cont})$, one has a  $(n-1)$-plane field $Q$ 
defined near $Z(M)$ with the following 
properties:
\begin{itemize}
\item at each point $z$ near $Z(M)$, the plane $Q_z$ is vertical, meaning that
 it is contained in the fibre of $TM$ passing through $z$;
\item $Q_z$ carries a symplectic bilinear form, well defined up to a positive factor;
\item as a (conformally) symplectic bundle over a neighborhood of $Z(M)$,  the plane field $Q$ 
is invariant by the holonomy of $\F$.
\end{itemize}
Then, there is a symplectic sub-bundle whose fibre at $x\in M$ is $S_x:= Q_{Z(x)}$ (the indeterminancy 
by a positive factor is irrelevant here). This is the underlying {\it formal contact structure}. 
 Another way to say the same thing consists of the following: $S$ is the kernel of a 1-form $\al$ and
  there is a 2-form $\beta$ on $M$ such that $\beta$ makes $S$ be a  symplectic  bundle; equivalently, it may be said that $\al\wedge\beta^{\frac{n-1}{2}}$ is a volume form.
 \item Assume $n=p+q$.
Given  $\F\in H^1_{\tau M}(M;\Gamma_{n,q}^{\rm fol})$, one has
 an
associated
 foliation $\G$ of codimension $q$ on a neighborhood of $Z(M)$
in the total space $TM$.
 The foliation 
  $\G$  induces on $M$ a $\Gamma_q$-structure $\gamma:=Z^*(\G)$,
  whose normal bundle is $\nu_\gamma:=Z^*(\nu_\G)$;
  and a monomorphism of vector bundles $\epsilon:\nu_\gamma\hookrightarrow \tau M$. 
  Indeed, for every point $x\in M$,
  the foliation $\G$ being transverse to $T_xM$ at $Z(x)$,
   the normal to $\G$ at $Z(x)$ embeds into $T_xM$. The pair $(\gamma,\epsilon)$,
   an \emph{augmented $\Gamma_q$-structure} {  (according to the vocabulary from \cite{wrinkling3}),} plays the role of a formal
   geometry associated to $\F$.
\item Assume $n\ge\dim Y$.
Given $\F\in H^1_{\tau M}(M,\Gamma_{n}^{Y})$, one has
an
associated submersion {   $w$ from a neighborhood of $Z(M)$
in the total space $TM$ to $Y$. This submersion
 $w$ induces a \emph{formal submersion} $(f,F)$ from $M$ to $Y$,
 that is, a bundle epimorphism from $TM$ to $TY$ whose value at every $x\in M$ is:
 $f(x)=w(Z(x))$, $F_x:=Dw_{Z(x)}\vert T_xM$.}
 \end{enumerate}
 Observe that all these spaces of formal geometries have natural topologies.
 Our first theorem yields a converse: for the above geometries, 
  formal $\Ga$-geometries lead
  to $\Ga$-structures.

 \begin{thm}\label{step_1_thm}
 {   Let $M$ be an $n$-dimensional manifold, possibly closed.  
 Let $\Ga$ be a groupoid in the set 
 of $n$-dimensional 
 geometries $\{\Gamma^{ \rm symp}_n,\Gamma_n^{\rm cont},\Gamma_{n,q}^{\rm fol},
\Gamma_{n}^{Y} \}$. } Then, the forgetful map  from 
 $H^1_{\tau M}(M; \Ga)$ to the corresponding space of formal $\Ga$-geometries  is a 
{     homotopy equivalence. }

  \end{thm}
  
  \begin{remarques} {\rm ${}$
 {   \nd 1) Actually, according to R. Palais   (\cite{palais} Theorem 15), the considered spaces have
  the property that a weak homotopy equivalence is a genuine homotopy equivalence. Thus, 
  it is sufficient to prove  that the mentioned forgetful map is a weak homotopy equivalence, meaning that 
  it induces an isomorphism of homotopy groups in each degree.\smallskip}
  
  \nd 2) In the case of symplectic/contact geometry, D. McDuff proved theorems of the same flavor
  using the {\it convex integration} technique of Gromov (\cite{mcduff}, 
see also \cite{elias} p. 104, 138).\smallskip

   \nd 3) Let $\F\in H^1_{\tau M}(M; \Ga)$.  { By taking} a section $s$ of $\tau M$ valued 
in the domain 
   foliated by $\F$ and generic with respect to $\F$, { there is an   induced  } $\Ga$-geometry 
   with singularities on $M\cong s(M)$.
   This seems to be a very natural notion of singular symplectic/contact structure. 
It follows from 
   Theorem \ref{step_2_thm} that the singular locus may be localized in a ball of $M$.
   }
  \end{remarques}

\begin{rien} {\bf Homotopy and regularization.} {\rm 
Our second theorem will allow us to regularize every parametric family of $\Gamma$-structures
 on every  manifold $M$ which is \emph{open},  
  that is, which has no closed connected component; this terminology
  will be permanently used in what follows.
 
 A \emph{homotopy} (also called a {\it concordance}\footnote{This second word emphasizes the difference with a one-parameter family of $\Gamma$-structures.}) 
 between two $\Gamma$-structures $(\nu_i,\F_i)$ ($i=0,1$)
on $M$ is a $\Gamma$-structure on $M\times[0,1]$
whose restriction to $M\times 0$ (resp. $M\times 1$) equals $(\nu_0,\F_0)$
(resp. $(\nu_1,\F_1)$). Of course, $\nu_0$ and $\nu_1$ must be isomorphic.

 A $\Gamma$-structure $(\nu,\F)$ is said to be 
 \emph{regular} if the foliation $\F$ is transverse not only to the fibres of $\nu$ but also 
 to $Z(M)$ in $E(\nu)$. This bi-transversality of $\F$
induces an isomorphism $\nu\cong\tau\bigl(Z(M)\bigr)$.
In that case, the pull-back $Z^*(\F)$ is 
a $\Gamma$-geometry on $M$, namely the foliation by points equipped with a 
transverse $\Ga$-geometry.

 }
\end{rien}

\begin{rien}{\bf The exponential $\Ga_n$-structure.} {\rm Given a complete Riemannian metric on 
the $n$-manifold $M$,
there is a well defined map 
$$exp: TM\to M.$$
When restricting $exp$ to a small neighborhood $U$ of $Z(M)$ in $TM$, we get a submersion to $M$. The foliation defined by the level sets of $exp\vert U$ represents a regular $\Ga_n$-structure on $M$,
denoted by $\F_{exp}\in H^1_{\tau M}(M; \Ga_n)$. Up to isomorphism (vertical isotopy in $TM$),
$\F_{exp}$ does not depend on the Riemannian metric as it is shown by the next construction.

Consider the product $M\times M$ and its diagonal $\De\cong M$. We have two projections
$p_v, p_h: M\times M\to \De$, respectively the vertical and the horizontal projection. A small tube $U$
about $\De$ equipped with $p_v$ is isomorphic to $\tau M$ as micro-bundle. Then,
 the same tube equipped with $p_h$ defines the $\Ga_n$-structure $\F_{exp}$\,.
}
\end{rien}

We recall the fundamental property of  the differential of $exp$ (independent of any Riemannian metric):
$$d (exp\vert T_xM)_{Z(x)}=Id : T_xM\to T_xM.
$$
As a consequence, if $f: M\to Y$ is a smooth map and $v\in T_xM$, one has 
\begin{equation}\label{diff}
f\circ exp_x(v)-f(x)= df_x(v) +o(\Vert v\Vert).
\end{equation}

\begin{thm}\label{step_2_thm}
Let $X$ be an $n$-manifold, let $\Gamma\subset\Gamma(X)$ be an open subgroupoid and 
let $M$ be a (connected) $n$-manifold.
Assume that $M$ is open (that is, no connected component is closed). Let {  
$$s\mapsto \xi_s=(\tau M, \F_s):\D^k\to H^1_{\tau M}(M;\Gamma)$$
 be a continuous family  of  tangential $\Ga$-structures, parametrized by the compact $k$-disk ($k\geq 0$),
  such that for every $s\in\partial\D^k$, the $\Gamma$-structure
 $\xi_s$ is regular and $\F_s$ is tangent to $\F_{exp}$ along $Z(M)$. 
  
 Then, there exists a continuous family { of concordances}
 $$s\mapsto\bar\xi_s=(\tau M\times[0,1], \bar\F_s):\D^k\to H^1_{\tau M}(M\times[0,1];\Gamma)$$
   such that
 \begin{itemize}
 \item 
 $\bar\F_s=pr_1^*(\F_s)$ for every $s\in\partial\D^k$, where $pr_1: M\times[0,1]\to M$ is the projection;
 \item  $\bar\F_s\vert(M\times 0)=\F_s$ for every $s\in\D^k$;
 \item for every $s\in\D^k$, the $\Gamma$-structure $\bar\xi_s\vert(M\times 1)$ 
 is regular and $\bar\F_s$ is tangent to $\F_{exp}$ along $Z(M\times 1)$.
 \end{itemize}
 }
\end{thm}

\begin{remarque}{\rm

The Smale-Hirsch classification of immersions $S\to Y$ (see \cite{smale,hirsch}), where $S$ is a closed manifold 
of dimension less than $\dim Y$, is covered by Theorem \ref{step_2_thm}; in particular, the famous {\it sphere
eversion} amounts to the case where $S$ is the 2-sphere, $Y= \R^3 $ and $k=1$. Let  us show it.

Let $(f,F):TS\to TY$ be a formal immersion. Then thanks to $F$ we have a monomorphism
$ F_*: \tau S\to f^*\tau Y$ over $Id_S$. Let $\nu$ be a 
complementary sub-bundle to the image of $ F_*$;  when $f$ is an immersion, $\nu$ is its normal bundle. 
Let $\hat S$ be a disk bundle in $\nu$; it is a compact manifold with non-empty boundary and 
$\dim \hat S=\dim Y$.
Thus, instead of immersing $S$ to $Y$ one tries to immerse $\hat S$ to $Y$; if it is done, the restriction to the 
0-section yields an immersion of $S$ to $Y$ with normal bundle $\nu$. The
 formal immersion $(f,F):TS\to TY$ easily
extends to a formal immersion $(\hat f, \hat F):\hat S\to Y$ in codimension 0.
Since $\hat F: T_x\hat S\to T_{\hat f(x)} Y$ is a linear isomorphism for every $x\in \hat S$, the level sets of 
$exp_Y\circ F$ is a $\Ga_n^Y$-foliation $\F$ near $Z(\hat S)$, that is, a $\Ga_n^Y$-structure on $\hat S$.
Moreover, thanks to Equation (\ref{diff}), if $\hat f$ is an immersion $\F$ is tangent to $\F_{exp}$; 
here $exp$ stands for  $exp_{\hat S}$.
Then, Theorem \ref{step_2_thm} applies and yields the desired immersion (or family of immersions).
}
\end{remarque}

\begin{cor}  Let $\Ga$ be a groupoid as in Theorem \ref{step_2_thm} and $\xi=(\tau M,\mathcal F)$ be a 
tangential $\Ga$-structure on a closed manifold $M$. Then, after a suitable concordance, all singularities
(that is,  the points where $\mathcal F$ is not transverse to $Z(M)\cong M$) are confined in a ball.
\end{cor}

\nd \proof
Let $B\subset M$ be a closed $n$-ball. Apply  Theorem \ref{step_2_thm}  to 
$\xi\vert (M\smallsetminus int\,B)$. We are given a regularization concordance $C$ of this 
restricted $\Ga$-structure. Since this concordance is given on a manifold with boundary, it extends 
to the whole manifold. Indeed, $B\times[0,1]$ collapses to 
$\bigl(B\times\{0\}\bigr)\cup\bigl( \partial B\times [0,1]\bigr)$.
\bull
{   \begin{remarque} {\rm Y. Eliashberg \& E. Murphy (\cite{murphy} Corollary 1.6) gave a similar result 
 for 
symplectic structures on closed almost symplectic manifolds of dimension greater than 4. Moreover,
  in the confining ball $B$ their singular symplectic structure is the {\it negative} 
cone of an {\it overtwisted } contact structure on $\partial B$.
Their proof  is based on the new techniques in contact geometry 
initiated by E. Murphy  \cite{murphy0} and developped in \cite{borman}.
}
\end{remarque}}

\begin{rien}{\bf The classical $h$-principle for $\Ga$-geometries.}
{\rm  For a groupoid $\Ga$ as  listed in 
\ref{list-geo} the $h$-principle states the following:}

If $M$ is an open $n$-manifold, 
the space of $\Ga$-geometries on $M$ has the 
same (weak) homotopy type as the space of formal $\Ga$-geometries on $M$.\smallskip

{\rm
This statement follows from Theorems \ref{step_1_thm} and \ref{step_2_thm}.\\

\nd \proof 
We start with a $k$-parameter family of formal $\Ga$-geometries on 
$M$, $k\geq 0$, which are genuine $\Ga$-geometries when the parameter  $s$ lies in $\partial \D^k$.
Then, for every $s\in \D^k$, the foliation $\F_{exp}$ is a $\Ga$-foliation near $Z(M)$. Thus, Theorem 
\ref{step_1_thm} applies and yields a $k$-parameter family of $\Ga$-structures on $M$ which remains 
unchanged when $s\in \partial\D^k$. Now, since $M$ is open, Theorem \ref{step_2_thm} applies and all the relative homotopy groups of the pair $(f\!ormal\ \Ga$-$geometries, \Ga$-$geometries)$ vanish.
\bull

}
\end{rien}

The article is organized as follows. In Section \ref{section_thom}, we detail the tool that goes back to 
R. Thom \cite{thom59}
and we prove Theorem  \ref{step_1_thm}
for submersion structures and for foliation structures.
 The next sections are devoted to the proof of Theorem \ref{step_1_thm} 
in the case of 
transversely symplectic structures. The existence part is treated in Section \ref{section_symp}. 
The family of such structures are considered in Section \ref{parametric_symp}; the proof of Theorem
\ref{step_1_thm}  is completed there when the groupoid is $\Ga_n^{ \rm symp}$. In Section 
\ref{section_contact},
 we adapt the proof to the groupoid $\Ga_n^{ \rm cont}$.
  Finally, in Section \ref{open}, we solve 
 the problem of regularizing the $\Ga$-structures on every open manifold.

\section{Thom's {   subdivision and 
jiggling}}\label{section_thom}

Reference  \cite{thom59} is the report of a lecture where R. Thom announced a sort of homological
$h$-principle (ten years before Gromov's thesis). 
A statement and a sketch of proof are given there; the details never appeared. From this text, we  extracted 
an  unusual subdivision process of the standard simplex and {   we derived two jiggling formulas\footnote{Thom speaks of 
 ``dents de scie'' 
 ({   saw teeth}); we keep the word {\it jiggling} that W. Thurston introduced in 
 \cite{thu74}.}. Our jiggling will be vertical 
 while Thom's jiggling is transverse to the fibres in some jet bundle. Nevertheless, we whall speak of 
 Thom's jiggling for it mainly relies on Thom's subdivision.} Actually, neither statement nor proof nor 
 formula were written 
 {   in  \cite{thom59}}, only words
 describing the object, a beautiful object indeed. 
 
 {   Here is a good occasion for mentioning that the famous {\it Holonomic Approximation Theorem}
 by Y. Eliashberg and N. Mishachev (\cite{elias} Chapter 3) is also based on a jiggling process, even 
 if that word is not used there.  The difference between their jiggling and ours is that the first one takes place in 
 the manifold itself while the second one is somehow vertical in the total space of a fibre bundle.
 }

 \begin{prop} \label{thom_sub}
 Let $\De^n$ denote the standard $n$-simplex. For every positive integer $n$, 
 there exist a non-trivial subdivision 
 $K_n$ of $\De^n$ and a simplicial map $\si_n: K_n\to \De^n$ such that:

\nd {\rm 1)} \emph{(non-degeneracy)}  the restriction of $\si_n$ to any $n$-simplex of $K_n$ is surjective;

\nd {\rm 2)} \emph{(heredity)} for any $(n-1)$-face $F$ of $\De^n$, 
 the intersection $K_n\cap F$ is { simplicially} isomorphic to $K_{n-1}$ 
 and $\si_n\vert F\cong \si_{n-1}$.
 \end{prop}
\proof Condition 2)  
 implies $\si_n(v)=v$ for any 
vertex of $\De_n$. For $K_1$, we may take $\De^1= [0,1]$ subdivised by two interior vertices: 
$0<v_1<v_0<1$ and we define $\si_1$ by $\si_1(v_1)=1$ and $ \si_1(v_0)=0$.

For $n=2$, let $A,B,C$ denote the vertices of $\De^2$. The polyhedron $K_2$ will be built
in the following way: subdivide each edge of $\De_2$ as $K_1$ subdivides $\De_1$;
add an interior triangle with vertices $a,b,c$ so that the line supporting $[b,c]$ is parallel to $[B,C]$
and separates  $A$ from $a$, etc.; join $a$ to 
the four vertices of $[B,C]$, etc. The simplicial map $\si_2$ is defined by 
$a\mapsto A, b\mapsto B, c\mapsto C$ and by imposing to coincide with $\si_1$ on each 
edge of $\De^2$. Condition 1)  is easily checked.

This construction extends to any dimension. If $K_{n-1}$ and $\si_{n-1}$ are known, each facet 
of $\De^n$ will be sudivided as $K_{n-1}$. Then, one puts a small $n$-simplex $\de^n$
in  the interior of $\De^n$ applying the same rules of parallelism and separation as for $n=2$.
Each vertex $v$ of $\de^n$ will be joined to the vertices of the facet $F(v)$ of $\De^n$
in front of $v$, 
this facet being sudivided by $K_{n-1}$. The map $\si_n$ 
maps $v$ to the vertex $V$ which is opposite to $F(v)$ in $\De^n$.\bull\\

\begin{remarques}\label{iteration}${}${\rm
\nd 1) This subdivision may be iterated $r$ times producing a subdivision $K_n^r$ which 
is arbitrarily fine and a simplicial map $\si_n^r: K_n^r\to \De^n$ fulfilling the two conditions 
of Proposition \ref{thom_sub}. More precisely,  thinking of $\si_n$ as a map from $\De^n$ to itself,
$\si_n^r$ will denote its $r$-th iterate and  $K_n^r$ is defined  by the next  formula: 
 \begin{equation}\label{folding}
        K_n^r= \left(\si_n^{r-1} \right)^{-1}(K_n)\,.  \end{equation}
We will call $\si_n^r$ an {\it $r$-folding map.}\smallskip

\nd 2) Thanks to 
 heredity {   (condition that  the barycentric subdivision does not fulfil)}, this subdivision of the standard simplex
and the $r$-folding map apply to any polyhedron.\smallskip

{   \nd 3) It is worth noticing that Thom's subdivision is not {\it crystalline} in the sense of H. Whitney
(\cite{whitney} Appencice II). Thus, it does not fit Thurston's techniques of {\it jiggling} (compare \cite{thu74}).
}}
\end{remarques}

Actually, the above construction has an unfolding property which is stated in the next proposition.

\begin{prop} \label{unfolding}
With the above notations, for every $n$-simplex $\tau$ of $K_n^{r-1}$, 
the restriction $\si_n^r\vert\tau$ is homotopic to $\si_n^{r-1}\vert\tau$ among piecewise linear 
 maps $\tau\to \De^n$ which are compatible with the face operators.
\end{prop}

\proof According to 
formula (\ref{folding}) 
it is sufficient  to prove the proposition
for $\si_n\equiv\si_n^1$. In that case, $\si_n^{r-1}= Id$. The homotopy is obvious for $n=1$; it consists of 
shrinking the middle interval $\de^1$ to the barycenter of $\De^1$ and shrinking its image
at the same time. Recursively, the homotopy of $\si_n$ is known on the faces of $\De^n$. Then, it is sufficient 
to define the homotopy on the interior small $n$-simplex $\de^n$. As when $n=1$,  the homotopy consists 
of shrinking $\de^n$ and its image simultaneously to the barycenter of $\De^n$. \bull

\begin{rien}{\bf First jiggling formula.} {\rm Let $M$ be an $n$-manifold and $\tau M=( TM
\mathop{\longrightarrow}\limits^p M)$ be its tangent bundle. Choose an auxiliary Riemannian metric
on $\tau M$ and an arbitrarily small open disk sub-bundle $U$ so that, for every $x\in M$,
the exponential map $exp_x: U_x\to M$ is an embedding. Take a combinatorial triangulation
$T$ of $M$ so fine that  every $n$-simplex $\tau$ of $T$ is covered by $exp_x(U_x)$ for every $x\in \tau$.
 Let $T^r$ be the $r$-th Thom subdivision of $T$ and $\si^r: T^r\to T$ be the 
 corresponding simplicial folding map. The $r$-th {\it jiggling map} $j^r:M\to TM$
is defined in the following way. For each $x\in M$, the point $j^r(x)$ is the unique point in $U_x$
such that 
\begin{equation}\label{1-formula} exp_x(j^r(x))=\si^r(x).
\end{equation}
This formula defines $j^r$ as a piecewise smooth section $M\to TM$.
We have the following properties.
}
\end{rien}
\begin{prop}${}$\label{jig}
\nd {\rm 1)} Let $\tau$ be an $n$-simplex of $T$, let  $x$ be a point in $\tau$ and let $\de$ be an $n$-simplex of $T^r$ passing through $x$.
Then, $j^r(\de)$ goes to $exp_x^{-1}(\tau)$ as $r\to\infty$. The convergence is uniform for $x\in \tau$.
 
 \nd {\rm 2)} The map $j^r$ is homotopic to the 0-section
$Z$ among  $PL$ maps which are transverse to the exponential 
foliation $\F_{exp}$ on each $n$-simplex of their domain. 

\end{prop}
\proof 
1) The diameter of the simplices of $T^r$ goes to 0 as $r$ goes to $\infty$. Then, for $y\in \de$, the point 
$j^r(y):= exp_y^{-1}\left(\si^r(y)\right)$ is close to $exp_x^{-1}\left(\si^r(y)\right)$. 
Since $\si^r$ is a surjective simplicial map onto $\tau$,
we have the $C^0$ closeness of $j^r(\de)$ and $exp_x^{-1}(\tau)$. A similar argument holds for the derivatives.\smallskip

\nd 2) 
On the one hand, the leaves of $\F_{exp}$ in  $U$ are $n$-disKs. We define $exp^u: U\to U$,
 $u\in [0,1]$, to be
 the map which is the homothety by $u$ in each fibre of $exp$. It is a homotopy from $Id_U$
 to $exp\vert U$ which restricts to a homotopy from $j^r$ to $\si^r$. On the other hand, 
 according to Proposition \ref{unfolding},
$\si^r$ is homotopic to $Id_M$ through $PL$ maps which are non-degenerate on each $n$-simplex 
of their domain,  
hence transverse to $\F_{exp}$.\bull

 \begin{remarques}\label{smoothness} ${}$
 
 \nd {\rm 1) Any piecewise smooth map defined on an $n$-manifold 
 $M$ and smooth on each $n$-simplex of a triangulation $T$  may be approximated 
 by a smooth map with the same polyhedral
 image. It is sufficient to precompose with a smooth homeomorphism such that, for every simplex
 $\tau$ in the $(n-1)$-skeleton of $T$ and every $x\in \tau$, 
  all partial derivatives in directions transverse $\tau$ {   vanish} at $x$.
  Then, even if  the concept of $\Ga$-structure is restricted to the smooth category, 
 there is no trouble to  pull-back a $\Ga$-structure by $j^r$;
  it will be well defined up to homotopy.\smallskip

  \nd 2) In general a jiggling, for instance based on the iterated barycentric triangulation,
  does not share the properties stated in Proposition \ref{jig} (non-degeneracy
and $PL$-homotopy).  }
 \end{remarques}
 
 \begin{rien}{\bf Second jiggling formula.} {\rm Here, we consider a trivial bundle $\ep^n$ of rank $n$ 
   whose base is an $n$-manifold $M$ equipped 
 with a {\it colored} triangulation\footnote{A triangulation $T$
 of dimension $n$
 is colored when each vertex has a color in $\{0,1,\ldots,n\}$ such that two vertices of the same simplex have different colors. The first barycentric subdivision of any {   triangulation} is colored. }.
 Let $\De^n\subset \R^n$ be a non-degenerate $n$-simplex whose vertices are colored. The coloring 
 defines a simplicial map $c:T\to \De^n$.
  We have a first jiggling  $j^1: M\to M\times\R^n$, $x\mapsto \bigl(x,c(x)\bigr)$. Then, the Thom process defines a $r$-th jiggling 
  \begin{equation}\label{2-formula}
  j^r(x)= \bigl(x, c\circ \si^r(x)\bigl).
  \end{equation}
 The first item of  Proposition \ref{jig} holds true for this formula: $j^r(\de)$ tends to $\{x\}\times\De^n$ when 
 $n$ goes to $\infty$. 
 }
 \end{rien}

\begin{rien}{\bf Proof of theorem \ref{step_1_thm} in the easy cases.}
{\rm
 For two of the four geometries considered, namely the submersion geometry
 $\Gamma^Y_n$ and the foliation geometry $\Gamma^{\rm fol}_{n,q}$, the jiggling
 method yields directly a simple proof of theorem \ref{step_1_thm}. 
 \medbreak
For the {\it submersion geometry},  we begin by proving that 
  the forgetful map is $\pi_0$-surjective.
 One is given  an $n$-manifold $M$, a $q$-manifold $Y$ 
 and a {\it formal submersion} (in the sense of  Subsection \ref{formal} (4)),
 that is, a pair $(f,F)$, where $f:M\to Y$ is {    a smooth} map and
 where $F:TM\to TY$ is a bundle epimorphism above $f$. 
 One seeks for a one-parameter family $(f_u,F_u)$ of formal
 submersions, $u\in[0,1]$, such that $(f_0,F_0)=(f,F)$, and such that $(f_1,F_1)$ underlies
  some $\Gamma_n^Y$-structure $\xi=(\tau M,\F)$. {   According the definitions given 
   in Subsections \ref{gamma_fol}  and \ref{formal},
  we have to find a pair $(w,\F)$ formed with  a submersion 
  valued in $Y$ and a codimension-$n$ foliation, both defined near $Z(M)$ in $TM$, such that:
  \begin{itemize}
  \item $w$ is constant on each leaf of $\F$;
  \item $f_1(x)= w(Z(x))$ and $F_1= Dw_{Z(x)}\vert T_xM$.
  \end{itemize}
  This will work by taking $w= exp\circ F$ which is clearly a submersion on some neighborhood 
  $U$ of $Z(M)$ in $TM$; here, 
  $Y$ is endowed with some auxiliary Riemannian metric and  $exp:TY\to Y$ is the associated exponential map. The only somehow delicate point is to find $\F$ as 
  a subfoliation of the foliation $\mathcal W$ whose leaves are $w^{-1}(y)$, $ y\in Y$.
  }
 Let $P$ be an $n$-dimensional plane field on $U$ transverse to every fibre $T_xM$ and contained in
 the kernel of the differential of {   $w$}.
 
  Let $T$ be a triangulation of $M$. Consider the iterated Thom subdivisions $T^r$.
   By Proposition
  \ref{jig} (1), for $r$ large enough, the $r$-th
   Thom jiggling $j^r$ maps every $n$-simplex of $T^r$ into $U$
and  transversely to $P$. Fix such an $r$.
 Then, on some small open neighborhood
  $V$ of $j^r(T^r)$ in $TM$, there is a $C^0$-small perturbation
  of $P$, among the $n$-plane fields  tangent to $\mathcal W$,
  yielding an integrable plane field on $V$. {   In the present situation where
   the dimension of the simplicial complex $j^r(T^r)$ is not larger 
  than the codimension of $P$, the wanted  integrating perturbation can be easily constructed 
  by induction on the dimension of the simplices (see the very beginning  of Section 6 in \cite{thu74}).
  }
  
     Let $S:M\to TM$
  be a smooth section so close to $j^r$, that $S(M)\subset V$. For every $u\in[0,1]$,
  one has the section {   $S_u:=uS$  valued in $U$}. 
 {   Set $f_u:=w\circ S_u$,
  and $F_u(v_x):=Dw_{uS(x)}v_x$. The structure  $S^*(w, \F)$ is really a $\Gamma_{n}^{Y}$-structure 
  whose underlying formal structure is $(f_1,F_1)$.}
  The $\pi_0$-surjectivity is proved.
  
  More generally, one is given a parametric family of formal submersions
   $(f_s,F_s)$, $s\in\D^k$, which  are underlying some
  $\Gamma_n^Y$-structures $\F_s\in H^1_{\tau M}(M;\Gamma_n^Y)$ for every 
  $s\in\partial\D^k$. First, one
  constructs, 
  on some open neighborhood $U$
  of the zero section in $TM$ and for every {   $s\in\D^k$, an $n$-plane field
  $P_s$ as above:
  $P_s$ is contained in $\ker D(exp\circ F_s)$,} and transverse to every fibre $T_xM$.
  One arranges that {   $P_s$ depends smoothly on $s$, and coincides with the tangent space 
  to $\F_s$ for every $s\in\partial
  \D^k$.}  The construction of such a family is easy, by convexity\footnote{   Compare footnote \ref{convex-note}.} of
  the space of the $n$-plane fields on $U$ contained in $\ker(D(exp\circ F_s))$ 
  and transverse to the fibres.
  
   Let $T$ be a triangulation of $M$, and let $T^r$ be a Thom 
subdivision whose order $r$ is large enough  so that the same jiggling $j^r(T^r)$
is transverse to {   $P_s$ for every $s\in \D^k$}; since the considered family is compact,
such an $r$ certainly exists. 
Then, the integrating
   perturbation can be chosen smoothly with respect to    $s$
  and coinciding with the identity for every $s\in\partial\D^k$. A single neighborhood $V$
  and a single section {   $S$ fit all  parameters $s\in\D^k$.}
  We get on $V$ a parametric family {   $(\PP_s)$ of $\Gamma_n^Y$-foliations
  transverse to every fibre $T_xM$; and $\PP_s=\F_s\vert V$
  for every $s\in\partial\D^k$.
  
  Define $S$ and  $S_u$ as above.
  Set $\F_s:=S^*(\PP_s)\in H^1_{\tau M}(M,\Gamma_{n}^{Y})$. 
Set    $f_{s,u}:=exp\circ F_s\circ S_u$
  and $F_{s,u}(v_x):=D(exp\circ F_s)_{uS(x)}v_x$. This is a one-parameter family
  of $\D^k$-parametrized families of formal submersions,
  between $(f_{s,0},F_{s,0})=(f_s,F_s)$ and $(f_{s,1},F_{s,1})$, which is
  the formal submersion underlying $\F_s$.

        Finally, the families $(\F_s)$ and $(\F_s)\vert\partial\D^k$
         are homotopic as mappings
        $\partial\D^k\to H^1_{\tau M}(M,\Gamma_{n}^{Y})$. The homotopy
        consists of  pulling  $\F_s$  back through $S_u$.  }
  The proof of Theorem \ref{step_1_thm} is complete for the groupoid $\Ga_n^{Y}$.\bull

In the case of the {\it foliation geometry} on a manifold $M$ of dimension $n=p+q$,
we are  given
 a parametric family of augmented $\Gamma_q$-structures
$(\xi_s,\epsilon_s)$, $s\in\D^k$. Moreover, for every {   $s\in\partial\D^k$, the augmented
$\Gamma_q$-structure $(\xi_s,\epsilon_s)$ is underlying some
  $\Gamma^{\rm fol}_{n,q}$-structure $\F_s\in H^1_{\tau M}(M;\Gamma^{\rm fol}_{n,q})$.
  
 Denote $\xi_s:=(\nu,
  \XX_s)$ this family of $\Gamma_q$-structures. Of course, the normal vector bundle $\nu$ over $M$ does not depend on
  $s\in\D^k$. Recall that $\epsilon_s:\nu\hookrightarrow TM$} is a monomorphism
  of vector bundles.
  
  For every  $s\in\partial\D^k$, denote by 
   $\G_s$  the foliation
  of codimension $q$ tangent to $\F_s$ on a neighborhood $U$ of $Z(M)$ in $TM${   ; 
  that is, if $\F_s$ is viewed as a codimension-$n$ foliation, we have $\F_s\subset \G_s$.}
  For every $x\in M$, define {   
  $\tau_s(x):=(T\G_s)_{Z(x)}\cap T_xM$ } to be the $p$-plane tangent
  to the foliation  {   $\G_s\cap T_xM$ at $Z(x)$.}
  Thus, {   $\tau M=\tau_s\oplus\epsilon_s(\nu)$.}
  The family $(\tau_s)$ extends to a $\D^k$-parametrized
  family $(\tau_s)$ of $p$-plane fields
  on $M$ complementary  to $\epsilon_s(\nu)$.
  
 For every $s\in\partial\D^k$,
 after pushing $\G_s$ by a vertical isotopy in $TM$, whose $1$-jet at every point of
  the zero section is the identity, and after restricting to some smaller neighborhood, 
  one can moreover assume  that
   the trace $\G_s\cap (T_xM\cap U)$ is the restriction to $T_xM\cap U$ of
   the linear $p$-dimensional foliation 
   parallel to $\tau_s(x)$. Then, the family
   $(\G_s)$ extends to a $\D^k$-parametrized family $(\G_s)$ of foliations of codimension $q$,
   transverse to the fibres of $TM\to M$ and 
   defined on some 
   neighborhood $U$ of $Z(M)$
  independent of $s$.
  Indeed,      for every $s\in\D^k
   \setminus\partial\D^k$,
   we define $\G_s$ as the pullback of $\XX_s$ through the linear projection of
   $TM=\tau_s\oplus\epsilon_s(\nu)$ onto $\epsilon_s(\nu)$ 
   parallel to $\tau_s$.
   
   Just as in the case of the submersion geometry, one constructs a
    $\D^k$-parametrized smooth family $(P_s)$
  of $n$-plane fields on $U$, transverse to the fibres $T_xM$,
   and contained in $T\G_s$. For $s\in\partial\D^k$, one has $P_s=T\F_s$.
   For a large enough integer $r$, a $C^0$-small perturbation
   of $P_s$ yields a foliation $\PP_s$ of codimension $n$,
   contained in $\G_s$, on some open
   neighborhood $V$ of the jiggled zero section $j^r(T^r)$.
   For $s\in\partial\D^k$, one has $\PP_s=\F_s$.
   
  Define  sections  $S$, $S_u$ as above.
   For every $s\in\D^k$, set $\F'_s:=S^*(\PP_s)
   \in H^1_{\tau M}(M,\Gamma_{n,q}^{\rm fol})$.
  The underlying augmented $\Gamma_q$-structure is homotopic to the given one
  $(\xi_s,\epsilon_s)$. The homotopy consists {   of} the
  $1$-parameter family of $\D^k$-parametrized families of augmented $\Gamma_q$-structures
  {   $(\nu,S_u^*(\G_s),
  \epsilon_s)$ (note that $\G_s$ is defined on the whole of the open set $U$).}
    Finally, the $\partial\D^k$-parametrized families $\F_s$ and $\F'_s$
    of $\Gamma_{n,q}^{\rm fol}$-structures are homotopic:
 the homotopy consists  {   of} $\F_{s,u}:=S_u^*(\F_s)$.
  The proof of Theorem \ref{step_1_thm} is complete for the groupoid $\Gamma_{n,q}^{\rm fol}$.\bull

}
\end{rien}

 \section{Existence of transversely symplectic $\Ga_n$-structures}\label{section_symp}
 
 In this section  we prove a slightly more general statement than the existence part  of Theorem \ref{step_1_thm} for 
 the groupoid $\Ga_n^{ \rm symp}$: we consider any  symplectic bundle of rank $n$. 
 We are going to use a more
 informative notation: a $\Ga_n^{ \rm symp}$-structure on $M$ will be denoted by $\xi=(\nu,\F,\Om)$
 where $\Om$ is a closed 2-form whose kernel is $\F$.
 
 \begin{thm}\label{existence}
  Let $\nu=(E\to M)$ be a symplectic bundle of  even rank $n$ over a manifold $M$ of dimension 
 $\leq n+1$.
 Then there exists a $\Ga_n^{ \rm symp}$-structure $\xi$
 on $M$ whose normal 
 bundle $\nu(\xi)$ is isomorphic to  $\nu$ as { a} symplectic bundle.

 Moreover, if  a real cohomology class $\bar a\in H^2(M,\R)$ is given, $\xi$ can be chosen 
 so that
 the cohomology class $[Z^*\Om]$ equals $\bar a$, where $\Om$ is the closed 2-form underlying  $\xi$.
 \end{thm}

 We think of this problem as a lifting problem that we attack by obstruction theory. Let us 
explain how it works. As for any groupoid of germs, 
there are a classifying space\footnote{The contravariant homotopy
functor $\Ga_n^{ \rm symp}(-)$ satisfies the axiom of gluing (Mayer-Vietoris) and wedge sum; then, the classifying space exists according to E. Brown's Theorem \cite{brown}.} $B\Ga_n^{ \rm symp}$ and 
 a canonical isomorphism $$\Ga_n^{ \rm symp}(M)\cong[M,B\Ga_n^{ \rm symp}] $$ 
where $[-,-]$ stands for the set of homotopy classes of maps.

 This classifying space is
 the source of two maps. The first one  is {   $\beta:B\Ga_n^{ \rm symp}\to B
Sp(n;\R) $:}
 if $f:M\to B\Ga_n^{ \rm symp}$ classifies a $\Ga_n^{\rm symp}$-structure $\xi=(\nu,\F,\Om)$
up to {concordance}, 
 $\beta\circ f$ classifies its normal bundle $\nu$. The second one is 
$\kappa: B\Ga_n^{ \rm symp}\to K(\R,2) $, where the target is the Eilenberg-MacLane space 
classifying the functor $H^2(-,\R)$: the composed map $\kappa\circ f$ classifies 
the cohomology class of  the closed 
2-form $Z^*\Om$.  
Finally, the {pair} 
$(\beta,\kappa)$ defines
a map 
$$\pi^{\rm symp}: B\Ga_n^{ \rm symp}\to BSp(n;\R)\times K(\R,2)
$$
{   that we see as  a homotopy fibration.} 
For Theorem \ref{existence}, we are given a map $M\to BSp(n;\R)\times K(\R,2)$
and we have to lift this map to $B\Ga_n^{ \rm symp}$.
Since $M$ is $(n+1)$-dimensional,  Theorem \ref{existence} is a direct corollary
 of the next statement (the $(n-1)$-connectedness would be sufficient for Theorem \ref{step_1_thm}).
{   Indeed, thanks to the long exact sequence associated with $\pi^{\rm symp}$, this map induces 
a monomorphism up to the $n$-th homotopy group (see, for instance Hatcher's book \cite{hatcher}, Section 4.3).}

\begin{thm} {\bf (Haefliger, McDuff)}\label{trans_symp} The homotopy fibre of $\pi^{\rm symp}$, denoted by $F\pi^{ \rm symp}$,
 is $n$-connected.
\end{thm}

 A. Haefliger (\cite{haefliger}, Section 6) showed that  the $(n-1)$-connectedness of this homotopy fibre is a consequence of the $h$-principle. 
 D. McDuff (\cite{mcduff}, Theorem 6.1) proved the $n$-connectedness thanks to the convex integration
 technique.

\begin{rien}{\bf  What do we have to prove for Theorem \ref{trans_symp}?}
{\rm We have to prove that the $k$-th homotopy group $\pi_k(F\pi^{ \rm symp})$ vanishes 
when $k\leq n$. An element of this group is represented by a  $\Ga_n^{ \rm symp}$-structure 
$\xi =(\ep^n,\F,\Om)$
on the $k$-sphere with the following properties:
\begin{itemize}
\item The normal bundle is trivial as a symplectic bundle; this means that its underlying 
symplectic bilinear form is the standard form {   $\om_0$} of $\R^n$ on each fibre. 
\item The {   closed} 2-form $\Om$, which is defined in a neighborhood of 
the 0-section $Z(S^k)$ in $\S^k\times\R^n$,
is assumed to be exact.
\end{itemize}

Let $(p_1,p_2)$ denote the two projections of $\S^k\times \R^n$ onto its factors.
Recall that the kernel of $\Om$ is the tangent space to the codimension-$n$
 foliation $\F$ and that $\F$ is transverse to the fibres of $p_1$.

We have   to 
extend this structure  $\xi$ over the $(k+1)$-ball $\D^{k+1}$ 
or, equivalently, to show that it is homotopic to the 
trivial structure $\xi_0:= (\ep^n,\F_0, \Om_0)$ where $\Om_0=p_2^*\,\om_0$.
 
According to Moser's Lemma with $k$ 
parameters \cite{moser}, 
 there exists a vertical isotopy  of $\S^k\times\R^n$, keeping $Z$ fixed, 
which reduces  us to the case 
where the germ at $Z(x)$ of the form induced by $\Om$ on $p_1^{-1}(x)$ equals 
 $\om_0$ for every $x\in\S^k$. After this {   vertical Moser isotopy,  
 take a  trivial tube $U=\S^k\times B^n$}
 in the domain of $\Om$, where $B^n $ is an $n$-ball of small radius. Now, Theorem \ref{trans_symp} directly follows from the next lemma,
  as we will see  just after its statement. 
}
\end{rien}
 \begin{lemme}\label{vertical} Given  the above-mentioned data, 
 there exist a section $s:\S^k\to U$, a neighborhood $W$ of $s(\S^k)$ in $U$
 and an ambient   diffeomorphism $\psi$ such that:
 \begin{enumerate}
\item $\psi$ is the time-$1$ map of a  \emph{vertical} isotopy $(\psi_t)$; set $W_0:= \psi(W)$;
\item $\psi$ sends the pair $(W, \Om)$ to $(W_0,\Om_0)$;
\item the isotopy $\psi_t$ is Hamiltonian with respect to $\om_0$ in each fibre.
 \end{enumerate}
 \end{lemme}
 Here, ``vertical'' means that the isotopy preserves each fibre of $p_1$.
 
 \begin{remarques}{ \rm ${}$
 
 \nd 1) The statement holds true  for every symplectic vector bundle of rank $n$,
  equipped (near the 0-section)
 with two forms exact forms $\Om$ and $\Om_0$ which define $\Ga_n^{\rm symp}$-foliations  and induce 
 the same symplectic  form on each fibre. 
 
 \nd 2) Moreover,
 the two first items  are valid for a pair of $\Ga_n$-foliations without any transverse geometry.}
 \end{remarques}

 \medskip
 \nd {\sc Proof of Theorem \ref{trans_symp}.} Since  $\psi$ is vertical, $s_0:= \psi\circ s$ is a section 
 of the trivial bundle $\ep^n$ over $\S^n$. Also, recall the 0-section $Z$.
 Then, we have a sequence of homotopies of $\Ga_n^{\rm symp}$-structures on $\S^k$: 
 \begin{itemize}
 \item a first homotopy from $Z^*\Om$ to $s^*\Om$;
 \item then, a homotopy from $s^*\Om$ to $s_0^*\Om_0$ defined by the isotopy $(\psi_t)$;
 \item a last homotopy from $s_0^*\Om_0$ to $Z^*\Om_0$.
 \end{itemize}
 The last structure obviously extends to the $(k+1)$-ball.\bull
 
{   Shortly said, the proof of Lemma \ref{vertical} will consist of taking a jiggled section in the sense of 
formula (\ref{2-formula})
whose simplices are very vertical, then 
covering it by boxes which {\it trivialize} the kernel of $\Om$ and pushing 
these boxes by some vertical Hamiltonian isotopy
until  $\ker\Om $ becomes horirontal. This isotopy is done
 recursively on the boxes. There are two main problems:
 \begin{enumerate}
 \item[-] rectifying the $(j+1)$-th box should not destroy what was gained for the $j$-th box;
 \item[-] manage the vertical isotopies to be Hamiltonian and not just symplectic; if not, they could 
  not extend.\\
 \end{enumerate}}
 
 \nd{\sc Proof of Lemma \ref{vertical}.} We limit ourselves to $k= n$; for $k<n$, it is the same argument 
 by replacing the base $\S^k$ with an $n$-dimensional base $\S^k\times \D^{n-k}$.  
Let $B^n:= p_2(U)$ and let $\De^n$ be a non-degenerate and 
 colored $n$-simplex in the interior of $B^n$.

 Take a decreasing sequence
 $$\ep_0>\cdots> \ep_j>\cdots >\ep_{n-1}.
 $$
 When $\al$ is a  strict closed {   $j$-face of $\De^n$, let $N(\al)$ denote 
the closed $\ep_j$}-neighborhood of $\al$ in $\R^n$. Set 
$$N(\De^n):= \De^n\mathop{\cup}\limits_\al N(\al)
$$where the union is taken over all faces of $\De^n$. For a suitable choice of the sequence $(\ep_j)$
we may arrange that :
 \begin{enumerate}
  \item  $N(\alpha)\cap N(\beta)=\emptyset$ if $\al$ and $\beta$ are two disjoint faces;
  \item if $\alpha\cap\beta\neq\emptyset$ and if $\alpha$ and $\beta$
  are not nested, then $N(\alpha)\cap N(\beta)$ is interior
  to $N(\alpha\cap\beta)$;
  \item $N(\De^n)\subset B^n$.
  \end{enumerate}
  Now, take a colored triangulation $T$ of the base $\S^n$, its Thom subdivision $T^r$ and the associated 
  jiggling $j^r$ given by formula (\ref{2-formula}).
  We are going to construct bi-foliated {\it boxes} associated with each simplex of 
  $T^r$ whose plaques are respectively contained in the leaves of $\F$ and in the fibres of $p_1$; 
  the boundary of a box has a part tangent to $\F$ and another part tangent to the fibres. 
  Let $\tau$ be a $k$-simplex of $T^r$; with $\tau$, the coloring of $T$ 
  associates some face $\tau^\perp$ of $\De^n\subset B^n$.
  The box $B(\tau)$ is defined  in the following  way. Its base $p_1(B(\tau))$
 equals  $star(\tau)$,
the star of $\tau$ in $T^r$.  In the fibre over  the barycenter $b(\tau)$, we take the domain $N(\tau^\perp)$. Finally, $B(\tau)$
is the union of all plaques of $\F$ passing through $N(\tau^\perp)$ and contained in $p_1^{-1}(star(\tau))$.
 If the diameter of the base 
is small enough, that is, if the order $r$ of the subdivision is large enough, the holonomy of $\F$
over the base is  $C^0$ close to {\it Identity}. 
Therefore, 
each plaque in $B(\tau)$ cannot get out of $U$; 
thus, it covers $star(\tau)$.

Look at two faces $\tau$ and $\tau'$ of the same simplex $\si$ of $T^r$. Assume first  that $\tau$ and $\tau'$
are disjoint. Apply  the above condition (1) to $\tau^\perp$ and $\tau'^\perp$; 
by the holonomy argument, if $r$ is large enough, the boxes $B(\tau)$ and $B(\tau')$ are disjoint.
Assume now that $\tau$ and $\tau'$ are not disjoint but not nested. Then, by (2), we have
$$B(\tau)\cap B(\tau')\subset B(\tau\cap \tau').
$$
Nevertheless, if $\tau$ and $\tau'$ do not belong to the same $n$-simplex and if
 $\partial star (\tau)\cap \partial star (\tau')\neq\emptyset$, then
$B(\tau)$ and $B(\tau')$ could intersect  badly. This is corrected  in the following way.

Again, for $r$ large enough, the leaves of $\F$ meeting $j^r(\tau)$ intersect the fibre over 
$b(\tau)$ in $N(\tau^\perp)$. This guarantees that  $j^r(T^r) $ is covered by the interior of the boxes.
From now on, $r$ is fixed. For  $1>\eta>0$, the {\it  $\eta$-reduced  box}
  associated with $\tau$ is defined by
$$B_\eta(\tau):= B(\tau)\cap p_1^{-1}\bigl((1-\eta)star(\tau)\bigr)$$
where the homothety  is applied from the barycenter $b(\tau)$. Fix  $\eta>0$ small enough
so that  the $\eta$-reduced open 
boxes
still cover the jiggling. Now, we are sure that $B_\eta(\tau)$ and $B_\eta(\tau')$ 
are disjoint once $\tau$ and $\tau'$ are disjoint. 

 The desired open set $W$ is the union $V_0\cup \cdots\cup V_k\cup \cdots\cup V_{n-1} $,
  where $V_k$ denotes the interior of the $\eta$-reduced  
 boxes associated with each $k$-simplex; the section $s$ is any smooth 
 approximation of  $j^r$ valued in $W$. We are ready to perform the isotopy. It is done step by step, 
 in the boxes associated with the vertices of $T^r$ first, then with the edges etc. For $x\in star(\tau)$,
 lifting the segment $[x,b(\tau)]$ to $\F$ yields 
 a holonomy  diffeomorphism  { between fibres of box}
 \begin{equation}\label{hol}
 (hol\,\F)_x^{b(\tau)}: B(\tau)_x\to B(\tau)_{b(\tau)}
 \end{equation}
which is an $\om_0$-symplectomorphism since $\Om$ is closed. Similarly, we have the
holonomy of $\F_0$ which also give an $\om_0$-symplectomorphism. The steps are numbered from 0 to 
$n$.

If $v$ is a vertex in $T^r$, we define $\psi^0$ in $B_\eta(v)$ by the next formula. For $z\in B_\eta(v)$
and $x=p_1(z)$, 
\begin{equation}\label{step1}
\psi^0(z)= (hol\,\F_0)_{v}^x\circ (hol\,\F)_x^{v}(z).
\end{equation}
Since the reduced boxes are disjoint, this formula simultaneously applies to  the reduced boxes associated 
with all vertices.  By shrinking the segment $[x, v]$ to $[x, x+t(v-x)] $ and by replacing $v$ with $x+t(v-x)$ in formula (\ref{step1}), we define an interpolation between $\psi^0(z)$ and $z$.
As a consequence $\psi^0$ is the time-1 map of a vertical isotopy
of embeddings $(\psi^0_t)$ 
which is easily checked to be symplectic. Since the components  of the domain of $(\psi^0_t) $ are contractible, 
this is actually a Hamiltonian isotopy\footnote{The infinitesimal generator $X_t$
of an $\om_0$-symplectic isotopy satisfies that  $\iota(X_t)\om_0$ is a closed 1-form; it is said to be Hamiltonian if this form is the differential of a function.}
which therefore extends to a global Hamiltonian isotopy
supported in $U$, still denoted by $\left(\psi^0_{t}\right)$.  Let $\F_1$ (resp. $\Om_1$) be the 
direct image of $\F$ (resp. $\Om$) by $\psi^0_{1}$; all reduced 
boxes are transported in this way, { becoming $B^1_\eta(\tau)$ for each $\tau\in T^r$.}
Observe that $\F_1$ is horizontal in the reduced new boxes  associated with  vertices.

The next step (numbered 1) deals with the edges. Let $e$ be an edge in $T^r$ with end points $v_0,v_1$. For 
$z\in B^1_\eta(e)$ and $x=p_1(z)$, define $\psi^1(z)$ by:
\begin{equation}\label{step2}
\psi^1(z)= (hol\,\F_0)_{b(e)}^x\circ (hol\,\F_1)_x^{b(e)}(z).
\end{equation}
 Observe that $\psi^1(z)=z$ when $z\in B^1_\eta(v_i),\ i= 0,1${; indeed, this box covers the barycenter $b(e)$ and $\F_1$ is horizontal there}. Moreover, $\psi^1$ is the time-1 map of
 a symplectic isotopy   $\left(\psi^1_{t}\right)$ 
 relative to the reduced boxes of the vertices; this isotopy, called 
 the {\it step-1 isotopy}, 
 follows from an  interpolation formula
 analogous to the one defining $\left(\psi^0_{t}\right)$. 
 
  If $e$ and $e'$ are two edges, after condition (2), the domain where their $\eta$-reduced boxes 
 could intersect
is contained in a domain where $ \F_1$ is horizontal and, hence, $\psi^1_{t}= Id$ on this domain. Therefore, 
$\psi^1_{t}$
is well defined on the union  $V_1$ of closed $\eta$-reduced boxes associated 
with the  vertices and  edges. Unfortunately, it is not a Hamiltonian isotopy of embeddings;
some vertical loops in $V_1$ may sweep out  some non-zero $\om_0$-area. 
Thus, it could  not extend  to an ambient vertical symplectic isotopy. {   The needed correction is offered by the next claim, following well-known ideas (compare V. Colin \cite{colin}, Lemme 4.4).}\\

\nd {\sc Claim.} {\it
{\rm 1)} There is a real combinatorial cocycle $\mu=\mu_{\Om_1}$ of the triangulation $T^r$  
such that, for each 
triangle $\tau$, the real number $<\mu,\tau>$ measures the $\om_0$-area swept out by the loop 
$\{x\}\times (\partial\tau^\perp)$  through the isotopy $\left(\psi^1_{t}\right)$ 
for every $x\in (1-\eta)star(\tau)$;
in particular, this area does not depend  on $x$.

 {\rm 2)} When $\Om_1$ is exact, $\mu$ is a coboundary.
 
{\rm 3)} There is an ambient vertical $\om_0$-symplectic isotopy $\left(g_t\right)_{t\in[0,1]}$,
 supported in $U$,
 which is stationary on $V_0$
and such that 
$\mu_{g_1^*\Om_1}=0$.
}
\medskip

The third item, together with the first item, means that  
the step-1 isotopy $\left(\psi^1_{t}\right)$ 
becomes Hamiltonian when $\F_1$ stands for the foliation tangent to 
$\ker g_1^*\Om_1$ instead of $\ker\Om_1$. 

The proof of the claim is postponed to the end of the section.
We first finish the proof of Lemma \ref{vertical} 
by applying the claim in the next way.

 After the step-0 isotopy, the cocycle $\mu_{\Om_1}$ is calculated and the Hamiltonian
 isotopy $(g_t)$ is derived.
Let $\tilde\F_1$ denote the foliation tangent to $\ker g_1^*\Om_1$; 
let $\tilde B^1_\eta(\tau):=\left(g_t\right)^{-1}\left(B^1_\eta(\tau)\right)$.
Now, the straightening formula \ref{step2} of the box $\tilde B^1_\eta(e)$ is applied
 $\tilde\F_1$ instead of $\F_1$. The associated isotopy $\left(\psi^1_t\right)$
 becomes  $\om_0$-Hamiltonian. Hence, it extends 
 to a vertical isotopy supported in $U$, denoted likewise,
  which is $\om_0$-Hamiltonian on each fibre $U_x$. This finishes 
  step 1 of the isotopy.
 
The next steps of the induction
are similar, except
that the question of being a Hamiltonian isotopy is not raised again since, up to homotopy,  every loop 
in $W$ is already contained in $V_0\cup V_1$. In the end of this induction, we have a proof of Lemma
 \ref{vertical} by taking $\psi=\psi^n_1$.
 \bull

 \bigskip
 
  \nd {\sc Proof of the claim.} 
 
 1) Let $e$ be an edge of $T^r$; its end points are denoted $v_0$ and $v_1$.
  Let  $x$ be a point in the base of $B_\eta(e)$.  Set $\ga:=\{x\}\times e^\perp$.
 We first compute the $\om_0$-area swept out by the vertical arc $\left(\psi^1_{1}\right)^{-1}(\ga)$  
 through the step-1 isotopy. Denote this area by $\mathcal A(x,e)$; 
 any other arc with the same end points would give the same  area.
 
 There are two  natural ``squares'', $C$ and $C_0$, appearing for this computation. The square $C$ 
 (resp. $C_0$)
 is generated by the holonomy of $\F_1$  (resp. $\F_0$) over $[b(e),x]$ with initial vertical arc $e^\perp$
 in the fibre $U_{b(e)}$. They have common horizontal edges: $\beta_i:=e\times v_i^\perp$ for $i=0,1$.
 Orient $e^\perp$ from $v_0^\perp$ to $v_1^\perp$; thus, $\ga$ and $\left(\psi^1_{1}\right)^{-1}(\ga)$ 
 are oriented
 by carrying the orientation of $e^\perp$ by the respective holonomies; and also 
 $C_0$ and $C$ are oriented by requiring $\{b(e)\}\times e^\perp$ to define the boundary orientation.
  Then, we have 
 \begin{equation}
 \mathcal A(x,e)= \int_C\Om_0-\int_{C_0}\Om_0\,.
 \end{equation}
 The second  summand is 0 by construction. Similarly, we have $\int_C\Om_1=0$. 
Then, if $\La$ is any primitive of 
 $\Om_1-\Om_0$, we derive 
 \begin{equation}\mathcal A(x,e)= -\int_C d\La\, .
 \end{equation}
 We now use a specific choice of  primitive.  Recall the zero-section $Z: \S^n\to U$. 
 For $t\in [0,1]$, let $c_t$
denote the contraction $(x,v)\mapsto (x, tv)$ and let $c: U\times [0,1]\to U$ 
be the corresponding homotopy
from $Z\circ p_1$ to $Id_U$.   This  yields 
the next formula:
\begin{equation}
\Om_1-\Om_0= d\left[ p_1^*\theta+ \int_0^1\! \iota\!\left({\partial_t}\right)c^*(\Om_1-\Om_0)\right]
\end{equation}
where $\theta$ is a primitive of the exact form $Z^*\Om_1$ (observe that $Z^* \Om_0=0$); 
the integral is just the mean value of a one-parameter family of 1-forms. This primitive 
of $\Om_1-\Om_0$  also reads 
\begin{equation}
\La_0:=p_1^* \theta +\int_0^1c_t^*\iota(v\partial_v)(\Om_1-\Om_0)\,,
\end{equation}
  which vanishes on every vertical vector  since $\Om_1$ and $\Om_0$ coincide on the fibres.
  Orient  $\beta_0$
 as the horizontal lift of $[b(e),x]$ and $\beta_1$ as the opposite of the oriented horizontal lift.
 We have
 \begin{equation}\int_C d\La_0= \int_{\beta_1}\La_0+\int_{\beta_0}\La_0\,.
 \end{equation}
 
 Now, we consider a triangle $\tau$ in $T^r$ and we look at the $\om_0$-area $\mathcal A(x,\partial\tau)$ 
 swept out by
  $\{x\}\times(\partial\tau)^\perp$ when $x$ belongs to $ (1-\de)star(\tau)$. 
 The vertices of $\tau$ are denoted by $v_i,\ i= 0, 1,2,$ {   cyclically} ordered;
  the oriented edges are $e_j:=[v_{j-1},v_j ]$ where $j-1$ is taken modulo 3. 
  There are two particular horizontal  lifts of $[b(e_{j}), x]$,  denoted by $\beta_{j,k}$ with $k=j$ or $j-1$ 
  depending on whether  its origin is $(b(e_{j}), v_j^\perp)$ or $(b(e_{j}), v_{j-1}^\perp)$. If $k=j-1$,
  it is oriented as $[b(e_{j}), x]$; if $k=j$, it has the opposite orientation.
   By summing up the area swept out by each edge of $\{x\}\times(\partial\tau)^\perp$, we have
  \begin{equation}
  \mathcal A(x,\partial\tau)=
  <\La_0, \beta_{1,1}+\beta_{2,1}+\beta_{2,2}+\beta_{0,2}+\beta_{0,0}+\beta_{1,0}>
  \end{equation}
 where the bracketing stands for the integration over chain. 
 
 Since $\Om_1-\Om_0$  
 vanishes on $(1-\eta)star(\tau)\times\{v_i^\perp\}$,  we have
 $$
 <\La_0, \beta_{i,i}+\beta_{i+1,i}>=<\La_0, [b(e_i),b(e_{i+1})]\times v_i^\perp>.
 $$
   By summation, we have
 \begin{equation}
  \mathcal A(x,\partial\tau)=\sum_i <\La_0, [b(e_i),b(e_{i+1})]\times v_i^\perp>
  \end{equation}
 which implies that $\mathcal A(x,\partial\tau)$ does not depend on $x$. 
 The combinatorial cochain $\mu$ is now defined by the next formula:
 \begin{equation}\label{cochain}
  <\mu,\tau>=\sum_i <\La_0, [b(e_i),b(e_{i+1})]\times v_i^\perp>.
  \end{equation}
  If an arbitrary primitive $\La$ of $\Om-\Om_0$ is used, the above formula
  becomes 
  \begin{equation}\label{cochain-2}
   <\mu,\tau>=\sum_i <\La, \{b(e_i)\}\times [v_{i-1}^\perp,v_i^\perp]> +<\La, [b(e_i),b(e_{i+1})]\times v_i^\perp>.
  \end{equation}
  Indeed, a change of primitive consists of adding a closed 1-form; and the integral of this on 
  the polygon $P$ considered in formula (\ref{cochain-2}) is zero since $P$ bounds a 2-cell.\footnote{The cochain $\mu$ is a cocycle.  Regarding the second item, this fact is not important
  and left to the 
  reader. Note that the previous calculation uses  a local primitive of $\Om$ only.}\\

  2) Since $T^r$ is a finite simplicial set, we only have to prove that   $<\mu,\Si>= 0$ for every 2-cycle 
   $\Si$ 
  of $T^r$. Here, the exactness of $\Om_1$ is used.
   Summing formula (\ref{cochain-2}) over all triangles of $\Si$ yields a sum of integrals 
  of $\La$ over horizontal polygons in regions where $d\La=0$ (one polygon for each vertex of $\Si$). 
 Then, these integrals are null. 
   Therefore, there exists a combinatorial 1-cochain $\al$ of $T^r$ such that $\mu=\partial^*\al$ where 
   $\partial^*$ stands for the  combinatorial co-differential. \\

  3) We are going to use this 1-cochain $\al$ in order to correct $\Om_1$ by a certain vertical isotopy. 
  Let $e$ be an oriented edge in $T^r$ with origin $v_-$ and extremity $v_+$. The value $\al(e)$ is used 
  in the following way. In the fibre over $b(e)$, we find an  $\om_0$-Hamiltonian isotopy 
   $(g_t^e)_{t\in [0,1]}$, compactly supported in $U_{b(e)}$ and 
   fixing  $(V_0)_{b(e)}$, 
  such that the area swept out by the arc $\ga_e:= b(e)\times [v_-^\perp,v_+^\perp]$
  is $-\al(e)$.\footnote{In dimension $n=2$, this is possible only if $\vert \la(e)\vert$ is less than
  the $\om_0$-area of $U_{b(e)}$. This last condition is satisfied when $r$ is large enough.}
   Observe that the Hamiltonian  function is not required to vanish in the fixed domain, but only to
  be constant on each connected component of the fiber $(V_0)_{b(e)}$ over $b(e)$.
  
 Then, the infinitesimal generator $X_t$ of the desired isotopy $(g_t)$ is chosen in finitely many fibres. By 
  a suitable partition of unity there is an extention which is Hamiltonian in each fibre, compactly 
  supported and vanishing in $V_0$. Note that the Hamiltonian has to be constant in each connected
  component of the fibre $(V_0)_{x}$, but these constants may vary with $x$.
  
  The 2-form $g_1^*\Om_1-\Om_1$ has a primitive associated with the isotopy, named
   the {\it Poincar\'e} primitive,
  \begin{equation}
  A= \int_0^1 g_t^*\left(\iota(X_t) \Om_1\right)dt.
  \end{equation}
 Since $X_t$ vanishes  on $V_0$, 
 the 1-form  $A$ vanishes over here and we have:
  \begin{equation}
  <A, \ga_e>= - \al(e).
  \end{equation}
  Now, $\La+A$ is a primitive of $g_1^*\Om_1-\Om_0$.  According to formula (\ref{cochain-2}),
  the combinatorial cochain $\mu_{g_1^*\Om_1}$ associated with the 2-form $g_1^*\Om_1$ vanishes
  and the claim is proved.
  \bull
 
\section{ Parametric family of transversely symplectic $\Ga_n$-structures }\label{parametric_symp} 

In this section, we prove the parametric version of Theorem \ref{step_1_thm}
for the groupoid $\Ga=\Ga_n^{\rm symp}$. We emphasize
that the required {   $k$}-connectedness of the homotopy fibre $F\pi^{\rm symp}$ depends only on the 
dimension of $M$ and not on the number of parameters in the family.  Indeed, there is no integrability 
condition with respect to the parameter. 
{   Moreover, we insist that a common jiggling will be used in the proof; its order is bounded by compactness of the parameter space.}

We consider the same setting as in Theorem \ref{existence}:   $\nu=(E\to M)$ is a bundle of
even rank $n$ over a manifold $M$ of dimension $\leq n+1$ equipped with a 
$k$-parameter family $(\om_u)_{u\in \D^k}$
of symplectic bilinear forms $\om_u$ on $E$. It is understood that $k$ is positive. 

\begin{thm} Assume there is a family $(\xi_u)_{u\in \partial\D^k}$
of $\Ga^{\rm symp}$-structures, namely a family  $(\Om_u)_{u\in \partial\D^k}$ 
of closed 2-forms defined near the zero section $Z$ of $E$, 
such that $\Om_u$ induces $\om_u$ on the fibres of $\nu$ for every $u\in \partial\D^k$\,\footnote{In other words, the symplectic normal bundles equal $(\nu,\om_u)_{u\in \partial\D^k}$. 
}.

Then,
this family extends over the whole $\D^k$ 
such that $\Om_u$ induces $\om_u$ on the fibres of $\nu$ for every $u\in\D^k$.
Moreover, the family of cohomology classes $[Z^*\Om_u]_{u\in \D^k}$ may be arbitrarily chosen among 
those which extend the boundary data.
\end{thm}

\proof We start with a cell decomposition  $\mathcal C$ of $M$ fine enough so that, 
for every $u\in \partial\D^k$
and every cell $C \in \mathcal C$, there is a fibered isotopy of $E_{\vert C}$
 (depending smoothly on $u$)  whose time-1 map $\psi_u$ satisfies:
$(\psi_u)_*\Om_u=\theta_u^*\Om_0$, where $\theta_u$ is a linear symplectic trivialization
of $(\nu_{\vert C},\om_u)$, depending smoothly on $u\in \D^k$, and where $\Om_0$ stands 
for the pull-back of $\om_0$ by the projection $C\times\R^n\to\R^n$.

The theorem will be proved by induction on an order of the simplices of $\mathcal C$ for which their 
dimension is 
a non-decreasing function. Skipping  the intermediate dimensions we jump to the $(n+1)$-cells.
Thus, we are reduced to consider the $n$-trivial bundle over $\S^n$ and a family 
$(\Om_u)$ of exact 2-forms on a small disK bundle $U$ about the zero section $Z$,
which induce the standard form $\om_0$ on each fibre $U_x$, $x\in\S^n$, (here a parametric version 
of Moser's lemma is applied again). This family fulfills the condition that $\Om_u=\Om_0$ for every 
$u\in \partial \D^k$. Let $T$ be a triangulation of $\S^n$ and let $T^r$ be a Thom 
subdivision whose order $r$ is large enough  so that the same jiggling $j^r(T^r)$
fits the proof of Lemma \ref{vertical} for every $u\in \D^k$; since the considered family is compact,
such an $r$ certainly exists. 

Each step of that proof may be performed with parameters using this fixed jiggling. Here it is worth
noticing that the vertical isotopy given by Lemma \ref{vertical} is stationary when $\Om_u=\Om_0$,
in particular when $u\in \partial \D^k$.\bull

The proof of Theorem \ref{step_1_thm} is now completed for the groupoid $\Ga_n^{\rm symp}$.\bull

\bigskip
\section{Transversely contact $\Ga_n$-structures}\label{section_contact}
Here, we prove a theorem which implies Theorem \ref{step_1_thm} for $\Ga_n^{\rm cont}$-structures.
{   Our setting is not the one of tangential $\Ga$-structures. It is the following.}
 Given an odd  natural integer $n$, a manifold $M$ and a vector bundle $\nu= (E\to M)$ of rank $n$, we 
 recall that  a  $\Ga_n^{\rm cont}$-structure on $M$ with normal bundle $\nu= (E\to M)$ is given by 
 $\xi=(A,\mathcal K)$, where  $A$ is a 1-form and $\mathcal K$ is a codimension-$n$ foliation,
  both  defined near the 0-section $Z$ in $E$, such that:
  \begin{itemize}
  \item $A\wedge dA^{\frac{n-1}{2}}$ induces a germ of volume form on $E_x$ for every $x\in M$;
  \item $\ker( A\wedge dA^{\frac{n-1}{2}})= T\mathcal K$;
  \item $\ker \!A$ contains $T\mathcal K$ and is invariant by the holonomy of $\mathcal K$\,.
  \end{itemize}
  
  {   As in the symplectic case, 
  the next Theorem was known to A. Haefliger \cite{haefliger} when $\dim M<n+1$ and to 
  D. McDuff \cite{mcduff} when $\dim M=n+1$.
  }

\begin{thm} \label{cont-structures}
Assume $M$ is a manifold of dimension not greater than the rank of  the vector bundle  $\nu$.  
Let $(\al,\beta)$ be  \emph{formal contact data}, that is,  a section $\al $ of $\nu^*$ and a section  $\beta$
of $\wedge^2\nu^*$ such that 
$\al\wedge \beta^{\frac{n-1}{2}}$ is  a non-vanishing section of {   $\wedge ^n \nu^*$.}
 Then, there exists a 
$\Ga_n^{\rm cont}$-structure $\xi= (A, \mathcal K)$ on $M$ with normal bundle $\nu$ such that,
for all $x\in M$, the next two conditions are fulfilled:
\begin{equation}\label{cont_equation} \left\{
\begin{array}{l}
 \ker \! A_{Z(x)}\cap \nu_x= \ker \al(x)\\
 (dA)_{Z(x)}= \beta(x).
\end{array}
\right.
\end{equation}

 Moreover, this statement holds true in a relative parametric version.
 \end{thm}
\proof For simplicity, we do not formulate any homotopy statement at the level of classifying spaces. Nevertheless, the strategy of  proof is similar to the one  we used for $\Ga_n^{\rm symp}$-structures. 
 It is  even simpler since every contact isotopy is Hamiltonian.

Let us first  consider  the non-parametric version. The construction of $\xi$ is performed step by 
step over each cell of a cell 
decomposition of $M$. We are looking on the last $n$-cell  $e^n$ only.
Let $C:= \partial e^n\times[0,1]\cong \S^{n-1}\times [0,1]$ be a collar neighborhood
of the boundary, on which we are given a $\Ga_n^{\rm cont}$-structure $\xi= (A,\mathcal K)$ which fufills
(\ref{cont_equation}).

Since the formal data $(\al,\beta) $ extends over $e^n$, there is a trivialization of $\nu$ over $e^n$,
$(p_1,p_2): E\vert e^n\to e^n\times \R^n$,
in which $(\al(x),\beta(x))$ is independent of $x\in e^n$. On $\R^n$ equipped with $(\alpha,\beta)$, 
we may think of $\beta $ as a closed differential form with constant coefficients; by taking a primitive,
we have a unique contact form $\al_0$ such that :
\begin{equation}\left\{
\begin{array}{l}
\al_0(0)=\al\\
d\al_0=\beta.
\end{array}
\right.
\end{equation}

We derive a {\it trivial} $\Ga_n^{\rm cont}$-structure $\xi_0=(A_0,\mathcal K_0)$ on $e_n$ 
such that
\begin{equation}\left\{
\begin{array}{l}
A_0=p_2^*(\al_0)\\
dA_0=p_2^*(d\al_0).
\end{array}
\right.
\end{equation}
Hence, $\mathcal K_0$ is the horizontal foliation.
Now, there is a Moser type lemma\footnote{The statement comes from Eliashberg-Mishachev's book
\cite{elias} where the proof is left to the reader. We only add the relative
and parametric version.
}
 which we are going  to present below.
This allows us to perform some vertical isotopy  which reduces to the case where, 
in a small tube $U$ about the zero section and for every point $x\in C$,
we have 
\begin{equation}\label{moser_cont-eq}
A_{\vert U_x}= {A_0}_{\vert U_x}
\end{equation}

\begin{lemme}\label{moser_cont}
 Let $\left(\al_t\right)_{ t\in [0,1]}$ be a path of contact forms in a manifold $V^{n}$.
Let $L$ be a hypersurface in $V$. It is assumed that the Reeb vector field $R_t$
of $\al_t$ is never tangent to $L$. Then, we have the following:

\begin{itemize}
\item[1)] The next equation whose unknown is $X_t$ can be solved near $L$:
\begin{equation}\label{5.1}
 L_{X_t}\al_t +\dot\al_t=0.
\end{equation} 

\item[2)] Let $\left(\al_t\right)_{ t\in [0,1]}$ be a path of germs in $(\R^{2p+1}, 0)$ of contact forms. If 
$\ker\al_t(0)$ is independent of $t$, then these germs are isotopic. 
\item[3)] The previous statements hold true with parameters and in a relative version.
\end{itemize}

\end{lemme}
\proof ${}$
\nd 1) The vector $X_t$  decomposes as $X_t=Y_t+Z_t$ with $Y_t\in \ker\al_t$
and 
$Z_t= \al_t(X_t)\, R_t$. Let us recall that $R_t$ generates in each point the 
kernel of $d\al_t$. Then, Equation (\ref{5.1}) becomes the following  system:
\begin{equation}\label{xx}
\left\{\begin{array}{c}
R_t\cdot\left(\al_t(X_t)\right)+\dot \al_t(R_t)=0\\
\iota(Y_t)\left({d\al_t}_{\vert\ker\al_t}\right)+ d \left(\al_t(X_t)\right)_{\vert \ker\al_t}+\dot{\al_t}_{
\vert\ker\al_t}=0
\end{array}
\right.
\end{equation}
Fix $t\in [0,1]$. The first equation of this system is a differential equation along the orbits of $R_t$
whose unknown function is $\al_t(X_t)$. It has a unique solution if $\al_t(X_t)$ is required  to equal  0
along $L$ (here the transversality assumption is used). Then, the component $Z_t$ of $X_t$ is determined.
 Once $\al_t(X_t)$ is known, the second equation of (\ref{xx}) has a unique solution since the form induced
 by $d\al_t$ on $\ker\al_t$
 is symplectic.\\

\nd 2) Here $L$ is the 
 hyperplane which is the common kernel of the contact forms in the considered path. Replace the germs
 with genuine representatives. Following the  solution of 1), we have
  $X_t(0)=0$ for every 
  $t\in [0,1]$. Therefore, the flow $\vp_t$ of $X_t$ keeps the origin fixed. It is well defined on some 
neighborhood of the origin 
up to $t=1$ and it is the identity on $L$. We deduce from 
Equation (\ref{5.1}) that the following is satisfied near the origin for every $t\in[0,1]$:
\begin{equation}\label{5.2}
\vp_t^*( L_{X_t}\al_t +\dot\al_t)=0,
\end{equation} 
and the latter is obtained by taking the time  derivative of  the next equation 
\begin{equation}\label{5.3}
\vp_t^*\al_t=\al_0\,.
\end{equation}
So, the desired isotopy is obtained 
by integrating $X_t$.  \\

\nd 3) Considering the equations which are solved, this claim is clear.\bull

We continue the proof of Theorem \ref{cont-structures}. In order to derive (\ref{moser_cont-eq})
from Lemma \ref{moser_cont}, we use $x\in C$ as a parameter and, in each fibre $E_x$, we consider
$\al_t= t A_{\vert E_x }+(1-t){A_0}_{\vert E_x}$. Due to the formal data, $\al_t$ is a contact form
near the origin of $E_x$ for every $t\in [0,1]$.

After this Moser type reduction, we have to state and prove a lemma  similar to
Lemma \ref{vertical}. Actually, it is not useful to write it down explicitly 
 since it is the same:
  the holonomy maps are contactomorphisms; thus, in each fibre  of a box the vertical isotopy preseves 
the contact distribution $\ker A_0\cap E_x$. Then, it is Hamiltonian with respect to $A_0\vert E_x$\footnote{
 The Hamiltonian function of a vertical isotopy of contactomorphisms whose infinitesimal generator is $X_t$
 is (in our setting) the time dependent function $z\in U \mapsto A_0(X_t)(z)$.
 }. 
Therefore, it extends globally since  extending such an isotopy 
amounts to extend its 
Hamiltonian function; thus, no obstruction is encountered. This finishes the proof of the non-parametric
version.

For the relative parametric version of Theorem \ref{cont-structures}, we have to
 consider a family $(\al_u,\beta_u)_{u\in \D^k}$
which  underlies a family of $\Ga_n^{\rm cont}$-structures when $u\in \partial\D^k$. Thanks to the relative parametric version of Lemma \ref{moser_cont},  we may  follow word for word the 
proof we gave for $\Ga_n^{\rm symp}$. 

${}$\bull

\section{Open manifolds}\label{open}

This section is devoted to the proof of Theorem \ref{step_2_thm}. Let us  recall the setting: $\Ga$ is an open subgroupoid of the structural groupoid $\Ga(X)$ 
of a model $n$-manifold $X$ (see Section \ref{intro});  
$M$ is an open manifold of dimension $n$ and  $\xi\in H^1_{\tau M}(M; \Ga)$ is a $\Ga$-structure 
on $M$ whose normal bundle is  $\tau M$, the tangent space to $M$. Let
{   $TM$ denote its total space and  $Z: M\to TM$ denote  the 0-section.}
The associated 
$\Ga$-foliation defined near $Z(M)$ in {   $TM$} is denoted by $\F=\F_\xi$. 
 We recall a topological fact about  open manifolds whose proof 
is available in {   \cite{elias} Section 4.3.} 

\begin{prop} \label{spine} 
Given an open $n$-manifold $M$, there exists an $(n-1)$-polyhedron $K\subset M$, called a
\emph{spine} of $M$, so that the inclusion is a homotopy equivalence. More precisely, for any 
 regular 
neighborhood $V$ of  $K$, there exists a \emph{compressing} isotopy  of embeddings $f_t: M\to M,\  t\in [0,1],$
from $Id_M$ to an embedding $f_1:M\to V$,  which is stationary on a neighborhood of $K$
and such that  $t'>t$ implies $f_{t'}(M)\subset f_t(M)$. 
\end{prop}

We are going to prove the next statement from which Theorem \ref{step_2_thm} will be 
easily
 derived.

\begin{thm}\label{spine_reg} Let $K\subset M$ be an $(n-1)$-dimensional polyhedron in an 
$n$-manifold $M$ (open or closed).
Let $\xi_s ,\  s\in \D^k$, be a $k$-parameter family 
 of $\Ga$-structures on $M$ with normal bundle $\tau M$. 
When $s\in \partial \D^k$, it is assumed that the associated foliation $\F_s$ is tangent to $\F_{exp}$
along $Z(M)$. Then, there exist an open neighborhood $V$ of $K$ in $M$ and a 
$k$-parameter family $\bar \xi_s$ 
of $\Ga$-structures on $M\times[0,1]$ 
-- that is, concordances of $\xi_s$ -- such that:
\begin{itemize}
\item $\bar\xi_s\vert V\times\{0\}= \xi_s\vert V$;
\item $\bar\xi_s\vert V\times\{1\}$ is regular  and its associated foliation is  tangent to $\F_{exp}$;
\item  $\bar\xi_s= p_1^*\left(\xi_s\vert V\right)$ for every $s\in \partial \D^k$, where $p_1$ denotes 
the projection $M\times[0,1]\to M$.
\end{itemize} 
\end{thm}
In other words, $(\bar\xi_s)_{s\in \D^k}$ is a family of  regularization concordances
on a neighborhood of a $K$,
 relative to the boundary of the parameter space.\\

\begin{rien}{\bf Proof of Theorem \ref{step_2_thm} from Theorem \ref{spine_reg}}. {\rm 
Here, $M$ is an open manifold. A spine $K$ of $M$ may be chosen 
(Proposition \ref{spine}) and 
Theorem \ref{spine_reg}  applies to these data: $K\subset M, (\xi_s)_{ s\in \D^k}$. 
So, we have a  family $\bar\xi_s$ of regularization 
concordances on some neighborhood $V$ of $K$ in $M$,
relative to $\partial \D^k$.
We have to extend this family to {   a family }
of regularization concordances 
 over the  whole of  $M$,
still relative to $\partial \D^k$. We may assume there exists $\rho$ close to 1 so that $\xi_s$
is regular on  $M$ when $\Vert s\Vert\in[\rho,1]$.

We first  insert the family of concordances 
described as follows,
where $t\in [0,1]$ is the parameter of the concordance: 
\begin{itemize}
\item for $\Vert s\Vert \leq \rho$, we put the concordance  $t\mapsto f_t^*\xi_s$, where $(f_t)_{t\in[0,1]}$
 is the isotopy of embeddings given by  Proposition \ref{spine};
\item for $\rho\leq\Vert s\Vert \leq 1$, we put the concordance
${ t \mapsto f_{\left(\frac{1-\Vert s\Vert}{1-\rho}t\right)}^*\xi_s}$.
\end{itemize}
{   When $t=1$ (that is the {\it end} of these concordances depending  on $s\in \D^k$)} and when $\Vert s\Vert\geq\rho$,
 the structures are regular on $M$. 
{   Then, denoting by $S: [1,2]\to [0,1]$ the shift $t\mapsto t-1$, we continue, for $t\in [1,2]$ with the concordances $(f_1\times S)^*\bar\xi_s$ when $\Vert s\Vert \leq\rho$; these ones are stationary
when $\Vert s\Vert =\rho$. Thus, we are   allowed } to extend them by the stationary concordances
when $\rho\leq \Vert s\Vert \leq 1$.  Of course, the previous piecewise description can be made smooth
if desired. \bull
}

\end{rien}
 \begin{rien} \label{claim}
 {\bf Proof of Theorem \ref{spine_reg} without parameters ($k=0$).}\label{no_parameter}
 
 {\rm  
 We start with a $\Ga$-structure $\xi$ on $M$. 
  Let   $ \F$ be its  associated $\Ga$-foliation 
 defined in some small neighbohood $U$ of $Z(M)$. Let $\xi^u$
 (resp.  $\F^u$) be the underlying $\Ga_n$-structure (resp. $\Ga_n$-foliation)
 of $\xi$ (resp. $\F$) where the transverse geometry is forgotten.

 The proof will consist of two steps: in the first step, we will make  a specific  
 regularization  
  of $\xi^u$ by some $\Ga_n$-concordance
 over $M\times [0,3]$; 
 in the second step, the geometric $\Ga$-structure of the concordance  will be defined 
 only over a small neighborhood of $K\times [0,3]$. 
  Finally, we get the $\Ga$-regularization of $\xi$ near $K$. \\
 
 \nd{\sc 1st step.}   Fix a small $\ep>0$. 
 As in Thurston \cite{thu74}, we consider a one-parameter family $P_t$, $t\in [0,3]$, of
  $n$-plane fields on $U$ with the following properties:
 
\begin{itemize}
 \item[--] $P_t$ is transverse to the fibres for every $t$. 
 \item[--]  $P_t$ is tangent to $\mathcal F$ when $t\in [0,1+\ep]$.
  \item[--] $P_t$ is  tangent to $\mathcal F_{exp}$ when $t\in [2, 3]$.
 \end{itemize}

  Such a plane field exists by barycentric combination  in the convex set\footnote{    \label{convex-note}Take an $n$-plane field $Q$ transverse to the fibres. The above-mentioned convex set is 
   affinely isomorphic to $\hom (Q,\tau^vTM)$
  where $\tau^v$ stands for the sub-bundle of $\tau(TM)$ tangent to the fibres of $TM\to M$.} of the 
  plane fields transverse to the fibres of $\tau M$.
   
   Let $T$ be a triangulation of $M$ containing a subdivision of  $K$ (also called $K$)
    as a sub-complex and 
    fine enough with respect to the open  covering 
   $\{exp_x(U_x) \mid x\in M\}$ in order that formula (\ref{1-formula}) makes sense. {   Here, we recall 
   that formula which holds for $x$ in any simplex of $T$:
   $$
   \exp_x(j^r(x)) =\si^r(x).
   $$
   }
   We now consider the Thom  jiggling given by  formula (\ref{1-formula}); 
   its order $r$   is chosen {   large enough} so that the
   $n$-simplices of $j^r(T^r)$ are transverse to 
  $P_t$ for every $t\in[0,3]$.

  The first piece of the  concordance, when $t\in [0,1+\ep]$, actually a $\Ga$-concordance,
  consists of moving the zero section from 
  $Z$ to $j^r$ by traversing any homotopy valued in $U$ and  
  stationary when $t\in[1,1+\ep]$. 
  The concordance of $\Ga$-structure is given 
  by pulling  $\xi$ back by this homotopy of maps $M\to U$
  (look at Remark \ref{smoothness}~1) about smoothness).

  We now describe the second piece of the concordance, when $t\in [1,2]$.
 We consider the  codimension $n$-plane field 
 $\tilde P$ in $U\times [0,3]$ defined by 
 \begin{equation}
 \tilde P(x,t): = P_t(x)\oplus \R \partial_t . 
 \end{equation}
 It is tangent to $\mathcal F\times [0,1+\ep]$ and to $\mathcal F_{exp}\times [2,3]$. 
   The trace of $\tilde P$ on each $(n+1)$-cell   of $j^r(M)\times [1,2]$
 is one-dimensional. Then, this trace is integrable. 
 Thus, there is a $C^0$-small smooth approximation of $\tilde P$, 
 relative to $t\in[0,1]\cup [2, 3]$ and still denoted
 by $\tilde P$, which is integrable near $j^r(M)\times [1,2]$. 
 Now,  the pair   $(j^r(M)\times [1,2],   \tilde P)$ 
 defines a concordance of $\Ga_n$-structure{   s}. This finishes the second piece. 
 
 The third  piece of the concordance when $ t\in[2,3]$ consists  of  keeping the foliation 
 $\F_{exp}$ and applying the homotopy
 from $j^r$ to the 0-section $Z$ as {   provided by} 
  Proposition  \ref{jig}~2). On the whole, we built  a specific regularization concordance
 of the underlying $\Ga_n$-structure $ \xi^u $, which is nearly sufficient  for our purpose.

 We need  more of {\it good position}. 
 Let $K^r   $ denote the  $(n-1)$-dimensional complex which is 
 the $r$-th Thom subdivision of $K$. 
  Let $\tilde K^r$  be the image of $Z(K^r)\times [0,3]$ along the  concordance built above. 
  This is an $n$-complex
 whose $n$-cells are not transverse to $\tilde P$. When $t\in [1,2]$,  the only reason 
 for non-transversality is 
 that $\tilde K^r$ and $\tilde P$ share the $\partial_t $-direction. 
 Let $\tilde K^r_{[t,t']}$ (resp. $\tilde K^r_t$) denote the  restriction 
 of $\tilde K^r$ over $M\times[t,t']$ (resp. $M\times\{t\})$.
 
 When  the $n$-cells  of $\tilde K^r_{[t,t']}$  are prismatic {   (that is, $\text{simplex}\times\text{interval}$)}
 which is always the case when $[t,t']\subset [1,2]$,
 they will receive the {\it standard subdivision} defined by H. Whitney (\cite{whitney}, Appendix II)\footnote{   This subdivision that  W. Thurston names {\it crystalline} is clearly explained inside the  proof of his famous {\it Jiggling} Lemma.}; this latter only depends on an order chosen on the set of vertices of $\tilde K^r_t$. \\

 \nd {\sc Claim.} {\it There exist a subdivision $t_1=1,t'_1,...,t_i, t'_i,..., t_N= 3$ 
 and a small piecewise smooth vertical isotopy, its time-one map being denoted by $\psi$, such that:
 \begin{itemize}
 \item[(i)] $\psi\vert \tilde K^r_{t_i}= Id$ for every $i= 1,...,N$;
 \item[(ii)] for every $n$-simplex  $\tau$ of the standard subdivision of $\tilde K^r _{[t_i,t'_i]}$ 
 (resp. $\tilde K^r _{[t'_i,t_{i+1}]}$), the image $\psi(\tau)$ is smoothly embedded in $U\times[0,1]$ and 
{   quasi-transverse}  to $\tilde P$. {   Here, quasi-transverse means transverse when $\dim\tau\geq n$
and no tangency when $\dim\tau<n$.}
 \end{itemize}
 }

  \nd {\sc Proof of the claim}. We search for a {\it jiggling in time}. We are going to do it 
  for $\tilde K^r _{[1,2]}$; the jiggling in time of $\tilde K^r _{[2,3]}$ is similar, but a bit more complicated 
  due to the fact that the cell decomposition is not purely prismatic. A numbering  of the vertices
  of $\tilde K^r _1$ is fixed: $v_1,v_2, \dots$; this numbering propagates to the corresponding 
   vertices of $\tilde K^r _t$ for every $t\in [1,2]$. 
  
    The time subdivision is chosen so that, for every $(n-1)$-simplex  
    $c\subset \tilde K^r_t$ and every  $x\in c$, 
  the hyperplane  $H_t(x,c):= T_x(c)+P_t(x)$
  varies very little in $TU$ when $t$ traverses  $[t_i,t_{i+1}]$, uniformly when  $x$ runs in any star.
   Let $t'_i$ be the middle of this interval. Each $\tilde K^r_{[t_i,t'_i]}$ (resp. $\tilde K^r_{[t'_i,t_{i+1}]}$)
  receives the {\it  standard subdivision} of the prismatic cells. On  $\tilde K^r_{[t_i,t'_i]}$, the desired embedding $\psi$ and
  its  isotopy
  from $Id$ are constructed recursively on the numbered stars of vertices $star(v_1), star(v_2),\dots$.
  Precisely, there is a locally finite family of isotopies $\chi_1, \chi_2,\dots $
  (that is, only finitely many supports intersect) and $\psi$ will be the composition of their time-one map:
  $\cdots\circ\chi_2^1\circ \chi_1^1$.
 The reversed isotopies are used over the interval $[t'_i,t_{i+1}]$.
  
  Let $v_1$ be the {\it first} vertex of  $\tilde K^r_{t_i}$. Let $X_1$ be a small vertical vector in $T_{v_1}U$
  which is chosen linearly independent  form all above-mentioned hyperplanes 
  $H_{t_i}(x,c)$ 
  where $c$ is any $n$-simplex of  $\tilde K^r_{t_i}$ passing through $v_1$.
  Let $v'_1$ be the corresponding vertex in  $\tilde K^r_{t'_i}$. By definition of the standard subdivision
  of $\tilde K^r_{[t_i,t'_i]}$, the vertex
  $v'_1$ is joined to $star(v_1, \tilde K^r_{t_i})$, the star of $v_1$ in $\tilde K^r_{t_i}$. In 
  affine notation $\chi_1^1$ is obtained by replacing in $\tilde K^r_{[t_i,t'_i]}$
  \begin{equation}
  star(v'_1, \tilde K^r_{[t_i,t'_i]})\quad \text{with}\quad \bigl(v'_1+X_1\bigr)*lk(v'_1,  \tilde K^r_{[t_i,t'_i]}),    \end{equation}
  where $lk(a,-)$ stands for the {\it link} of $a$, that is, the union of simplices in 
  $star(a,-)$ which do not contain $a$.   
  After this step, we get the property that every $n$-simplex in 
  $B_1:= \chi_1^1\left(\tilde K^r_{[t_i,t'_i]}\right)$ passing through $v'_1+X_1$ is transversal to $\tilde P$.
  
  The second step consists of a similar construction related to the star of 
    $v'_2$ in $B_1$ by using a small vertical vector $X_2$ which is linearly independent of the 
    hyperplanes
  made with the $n$-simplices of $B_1$ passing through $v'_2$. This yields $\chi_2^1$, the time-one map
  of the second isotopy which allows us to gain the transversality of new $n$-simplices to $\tilde P$.
  When $star(v'_2,B_1)$ meets
  $star(v'_1+X_1, B_1)$, the vector $X_2$ is chosen so small that the property gained in the first step 
  is preserved. And so on.
  \bull
  
    In what follows, we still denote  by $\tilde K^r$ the outcome of the previous jiggling. As a result, 
  every $n$-cell of $\tilde K^r_{[1, 3]}$  is transverse to $\tilde P$.
  This simplicial complex $\tilde K^r$ 
   collapses\footnote{A simplicial complex $L$
  collapses to $K$ if there is a sequence of elementary collapses $L_q\searrow L_{q+1}$
 starting with $L$ and ending with $K$. An elementary collapse means that $L_q$ is the union of 
 $L_{q+1}$ and a simplex $\si$ so that $\si\cap L_{q+1}$ is made of the boundary of $\si$ with an
 open facet removed. } 
 successively  to  $\tilde K^r_{[0,2]}$ and then to $\tilde K^r_{[0,1]}$. }\\
 \end{rien}
 \vfill\eject
 
 \nd{\sc 2nd step.} 
 
 We now focus on $\tilde K^r_{[0,1]}$ on which $\tilde K^r$ collapses.
     Since the cells of $\tilde K^r_{[1,3]}$ of positive dimension are quasi-transverse to the foliation $\tilde P$, 
  the collapse  $\tilde K^r \searrow \tilde K^r_{[0,1]}$
 extends to a collapse of pair
 \begin{equation}\label{collapse}
 (\tilde K^r, \tilde P)\searrow(\tilde K^r_{[0,1]},  \tilde P)\,.
 \end{equation}
 Let $N(\tilde K^r)$ denote a small neighborhood of $\tilde K^r$ in $U\times [0,3]$.
 From the sequence of elementary collapses, one derives step by step
 an embedding of pair
 \begin{equation}\label{emb}
 \Phi: \bigl(N(\tilde K^r), \tilde P\bigr)\to \bigl(U\times [0,1+\ep), \,\mathcal F\times [0,1+ \ep)\bigr)
 \end{equation}
 which induces  the inclusion $N(\tilde K^r_{[0,1]}) \hookrightarrow U\times[0,1+\ep)$.
 Since $\mathcal F\times [0,1+\ep]$ 
 is a $\Ga$-foliation, $\tilde P_{\vert N(\tilde K^r)}$ is so
  by pulling back through $\Phi$. Therefore, 
 $\Phi^*\bigl(\mathcal F\times[0,1+\ep)\bigr)$ is a regularization concordance 
  of the $\Ga$-structure which is induced  near   $K$. 
 \bull \\
 This last process associated with collapses is named {\it inflation} in W. Thurston's article \cite{thu74}.
 \begin{rien}{\bf Relative parametric version of Theorem \ref{spine_reg}.} {\rm
 Here, the data consist of a family $(\xi_s)_{s\in\D^k}$ of $\Ga$-structures whose normal bundle is
 the tangent bundle $\tau M$. It is assumed that,
 for every $s\in \partial \D^k$, the associated foliation $\F_s$ is tangent to $\F_{exp}$
 along the 0-section $Z(M)$, 
 hence  
 $\xi_s$ is regular. Without loss of generality, we may assume $\xi_s$ is regular 
 when $1\geq \Vert s\Vert\geq \rho $ for some $\rho$ close to 1.
 The proof just  consists of two remarks.
 
 1) The previous proof (see Subsection \ref{no_parameter})
  works directly for our $k$-parameter family of data  if we do not take care of the boundary condition. 
  Indeed, 
  observe that a common order {   $r$ of Thom jiggling $j^r$ can}
be chosen for all $s\in \D^k$ since the family of $n$-plane fields $P_{s,t}$ we have to look at
 is compact. {   Thus}, if the jiggling is vertical enough{  , that is, for $r$ large enough,} its $n$-simplices are transverse to $P_{s,t}$
 for every $(s,t)\in \D^k\times [0,3]$. The 0-parameter process applies for every $s\in \D^k$
 and yields  a regularization in a fixed neighborhood  $V$ of $K$ in $M$. Precisely, we have formulas
  (\ref{collapse}) and (\ref{emb}) depending on the parameter $s\in \D^k$, yielding 
  regularization concordances $\Phi_s^*\bigl(\mathcal F_s\times[0,1+\ep)\bigr)$.

 2)  We may assume  that $P_{s,t}$ is tangent to $\F_s$ for every 
 $s\in\{1\geq  \Vert s\Vert\geq \rho\}$ and $ t\in [0,3]$.  For $\Vert s\Vert\in [\rho,1]$,  
  set $\mu(s):= \frac{1- \Vert s\Vert}{1-\rho}$. Recall that  $\Phi_s$ is the time-one map of an isotopy 
  of embeddings $\Phi_s^w: N(\tilde K^r)\to U\times [0,3] $, $w\in [0,1]$, relative to 
  $N(\tilde K^r_{[0,1]})$ and such that  
  $\Phi_s^0= Id$ and $\Phi_s^1\bigl(N(\tilde K^r)\bigl) \subset N(\tilde K^r_{[0,1+\ep]})$.
   
   We finish, for $s\in\{1\geq  \Vert s\Vert\geq \rho\}$, with the regularization concordance
  $(\Phi_s^{\mu(s)})^*\bigl(\mathcal F_s\times[0,3]\bigr)$. When $\Vert s\Vert =1$, this is the {trivial }
  concordance $\mathcal F_s\times [0,3]$. Then, the relative version is proved. 
   \bull
  }
 \end{rien}
 
 {   \nd{\bf Acknoledgments.} We thank Howard Jacobowitz and Peter Landweber who have expressed
  their interest towards our article and sent to us some valuable remarks.
  We are deeply grateful to the anonymous referee for his exceptionally careful reading and
  for having suggested to us  many   improvements of writing.}


\vskip 1cm


\begin{thebibliography}{999}
{   \bibitem{borman} M. S. Borman, Y. Eliashberg \& E. Murphy, {\it Existence and classification of overtwisted contact structures in all dimensions,} arXiv:1404.6157.}
\bibitem{brown} E. Brown, {\it Cohomology theories,} Annals of Math. 75 (1962), 467-484.
{   \bibitem{colin} V. Colin, {\it Livres ouverts en g\'eom\'etrie de contact,} in {\it S\'eminaire  Bourbaki, 
vol.  2006/2007,} Ast\'erisque 317, Soc. Math. France (2008),  91-117.}
{  
\bibitem{wrinkling3} Y. Eliashberg \& N. Mishachev,{\it Wrinkling of smooth mappings III: Foliations of codimension greater
than one,} Topol. Methods Nonlinear Anal. 11 n$^\circ$2 (1998), 321-350.
\bibitem{elias} Y. Eliashberg \& N. Mishachev, {\it Introduction to the $h$-principle}, G.S.M 48, Amer. math. Soc., 2002.

 \bibitem{murphy}  Y. Eliashberg \& E. Murphy, {\it Making cobordisms symplectic,} arXiv:1504.06312.}
\bibitem{haefliger70} A. Haefliger, {\it Feuilletages sur les vari\'et\'es ouvertes}, Topology 9 (1970), 183-194.
\bibitem{haefliger} A. Haefliger, {\it Homotopy and integrability}, 133-175 in:
Manifolds-Amsterdam 1970, L.N.M. 197, Springer, 1971.
\bibitem{haefliger2} A. Haefliger, private communication, Jan. 2014.

\bibitem{gromov} M. Gromov, {\it Partial Differential Relations,} Springer, 2002.
1980. 
{  \bibitem{hatcher} A. Hatcher, {\it Algebraic topology,} Cambridge University Press, 2002. }
 \bibitem{hirsch} M. Hirsch, {\it Immersions of manifolds,} Trans. Amer. Math. Soc. 93 (1959), 242-276.
 \bibitem{mcduff} D. McDuff, {\it Application of convex integration to symplectic and contact geometry,}
 Ann. Inst. Fourier (Grenoble) 37  (1987), 107-133.
\bibitem{moser} J. Moser, {\it On the volume elements on a manifold,} Trans. Amer. Math. Soc.
120 (1965), 286-294.
{   \bibitem{murphy0} E. Murphy, {\it  Loose Legendrian embeddings in high dimensional contact manifolds, arXiv:1201.2245.}}
\bibitem{palais} R. Palais, {\it Homotopy theory of infinite dimensional manifolds}, Topology 5 (1966),
1-16.
\bibitem{smale} S. Smale, {\it The classification of immersions of spheres in euclidean spaces}, Ann.
of Math. 69 (1959), 327-344.

\bibitem{thom59} R. Thom, {\it Remarques sur les probl\`emes comportant des in\'equations diff\'erentielles
globales,} Bull. Soc. Math. France 87 (1959), 455-461.
\bibitem{thu74} W. Thurston, {\it The theory of foliations of codimension 
greater than one}, Comment. Math. Helv. 49 (1974), 214-231.
\bibitem{veblen_31} O. Veblen \&J. H. C. Whitehead, {\it A set of axioms for differential geometry,} Proceedings of the National Academy of Sciences 17, 10 (1931), 551-561.
\bibitem{whitney}   H. Whitney, {\it Geometric Integration Theory,} Princeton Univ. Press, 1957.
\end{thebibliography}
\end{document}